\documentclass[11pt]{article}
\usepackage{latexsym}
\usepackage{amsfonts,amssymb,amsmath}
\usepackage{amscd}
\newtheorem{thm}{Theorem}[section]
\newtheorem{cor}[thm]{Corollary}
\newtheorem{lem}[thm]{Lemma}
\newtheorem{pro}[thm]{Proposition}
\newtheorem{defn}[thm]{Definition}

\title{The Thompson-Higman monoids $M_{k,i}$ : 
the $\cal J$-order, the $\cal D$-relation, and their complexity}

\author{ Jean-Camille Birget  }

\date{\today}
\begin{document}
\maketitle

\begin{abstract}
The Thompson-Higman groups $G_{k,i}$ have a natural generalization to
monoids, called $M_{k,i}$, and inverse monoids, called ${\it Inv}_{k,i}$. 
We study some structural features of $M_{k,i}$ and ${\it Inv}_{k,i}$ and
investigate the computational complexity of related decision problems. 
The main interest of these monoids is their close connection  with 
circuits and circuit complexity. 

The maximal subgroups of $M_{k,1}$ and ${\it Inv}_{k,1}$ are isomorphic
to the groups $G_{k,j}$ ($1 \leq j \leq k-1$); so we rediscover all the
Thompson-Higman groups within $M_{k,1}$.

Deciding the Green relations $\leq_{\cal J}$ and $\equiv_{\cal D}$ of
$M_{k,1}$, when the inputs are words over a finite generating set of
$M_{k,1}$, is in {\sf P}. 

When a circuit-like generating set is used for $M_{k,1}$ then 
deciding $\leq_{\cal J}$ is {\sf coDP}-complete (where {\sf DP} is the 
complexity class consisting of differences of sets in {\sf NP}). 
The multiplier search problem for $\leq_{\cal J}$ is 
{\sf xNPsearch}-complete, whereas the multiplier search problems of 
$\leq_{\cal R}$ and $\leq_{\cal L}$ are not in {\sf xNPsearch} unless 
{\sf NP} $=$ {\sf coNP}. We introduce the class of search problems
{\sf xNPsearch} as a slight generalization of {\sf NPsearch}.

Deciding $\equiv_{\cal D}$ for $M_{k,1}$ when the inputs are words over a
circuit-like generating set, is 
$\oplus_{k-1} \! \bullet \! {\sf NP}$-complete; 
for any $h \geq 2$, $\oplus_h \! \bullet \! {\sf NP}$ is a modular counting 
complexity class, whose verification problems are in {\sf NP}.
Related problems for partial circuits are the image size problem (which is 
$\# \bullet {\sf NP}$-complete), and the image size modulo $h$ problem
(which is $\oplus_h \! \bullet \! {\sf NP}$-complete). 
For ${\it Inv}_{k,1}$ over a circuit-like generating set, deciding 
$\equiv_{\cal D}$ is $\oplus_{k-1} {\sf P}$-complete.
It is interesting that the little known complexity classes {\sf coDP}
and $\oplus_{k-1} \! \bullet \! {\sf NP}$ play a central role in $M_{k,1}$.
\end{abstract}

%%%%%%%%%%%%%%%%%%%%%%%%%%%%%%%%%%%%%%%%%%%%%%%%%%%%%%%%
% Section 1
%%%%%%%%%%%%%%%%%%%%%%%%%%%%%%%%%%%%%%%%%%%%%%%%%%%%%%%%

\section{Introduction}

The Thompson-Higman groups $G_{k,i}$, introduced by Graham Higman in
\cite{Hig74}, can be generalized in a straightforward way to monoids, 
denoted $M_{k,i}$, and inverse monoids, denoted ${\it Inv}_{k,i}$. 
The generalization of $G_{k,1}$ to $M_{k,1}$ and ${\it Inv}_{k,1}$, was
given in \cite{BiThomMon}. The definition of $M_{k,i}$ for $i \geq 2$ is 
a straightforward combination of the definitions of $M_{k,1}$ and 
$G_{k,i}$.  In brief, $M_{k,i}$ consists of all maximally extended right 
ideal homomorphisms between right ideals of $BA^*$, where $A$ and $B$ are 
finite alphabets with $|A| = k \geq 2$ and $|B| = i \geq 1$. Detailed 
definitions of $M_{k,i}$ and ${\it Inv}_{k,i}$ appear below.

This paper is a continuation of our study of monoid generalizations of the
Thompson-Higman groups. As in \cite{BiThomMon, BiRL}, our motivations are 
the following: \ (1) The generalization of $G_{k,i}$ to a monoid or an 
inverse monoid is natural and straightforward; \ (2) the monoids $M_{k,i}$ 
and ${\it Inv}_{k,i}$ have interesting and surprising properties; 
 \ (3) for certain infinite generating sets, the elements of $M_{2,1}$ are 
similar to circuits, with word-length polynomially equivalent to 
circuit-size.   

\smallskip

The definition of $M_{k,i}$ requires some preliminary notions, most of 
which are familiar from formal language theory, or information theory, or
algebra.
Let $A$ and $B$ be finite alphabets with $|A| = k \geq 2$ and 
$|B| = i \geq 1$. By $A^*$ we denote the set of all words over $A$, 
including the empty word $\varepsilon$. A {\it right ideal} of $A^*$ is 
any set $R \subseteq A^*$ such that $R = RA^*$. 

We consider the set $BA^*$, i.e., the set of all words of the form $b_jx$ 
with $b_j \in B$ and $x \in A^*$. Although $BA^*$ is not a monoid with 
respect to concatenation, we can nevertheless define the concept of a 
{\it right ideal} of $BA^*$: It is any set of the form $B_0R$, where 
$B_0 \subseteq B$ and where $R \subseteq A^*$ is any right ideal of $A^*$. 
A right ideal $R$ is {\it essential} iff all right ideals of $BA^*$ 
intersect $R$. (We say that two sets $S_1$ and $S_2$ {\it intersect} iff  
$S_1 \cap S_2 \neq \varnothing$.)
For right ideals $R_2 \subseteq R_1 \subseteq BA^*$, we say that $R_2$ is 
{\it essential in} $R_1$ iff all the right ideals that intersect $R_2$ also 
intersect $R_1$. Two right ideals $R_2$ and $R_1$ of $BA^*$ are 
{\it essentially equal} iff every right ideal of $BA^*$ that intersects
$R_2$ intersects $R_1$, and vice versa; in that case we write 
 \ $R_2 \ =_{\sf ess} \ R_1$. If $R_2 =_{\sf ess} R_1$ then
 \ $R_2 =_{\sf ess} R_1 =_{\sf ess} R_1 \cap R_2$.

A {\it prefix code} in $BA^*$ is any set $P \subseteq BA^*$ such that no 
word in $P$ is a prefix of another word in $P$; hence, a prefix code of 
$BA^*$ is of the form $B_0 Q$ for some $B_0 \subseteq B$ and some prefix 
code $Q \subseteq A^*$. 
A set $P \subset BA^*$ is  a {\it maximal prefix code} iff $P$ is a prefix 
code which is not a strict subset of any other prefix code in $BA^*$.

A {\it right ideal homomorphism} over $BA^*$ is a total surjective function 
$\varphi: R_1 \to R_2$ such that $R_1, R_2$ are right ideals of $BA^*$, and 
such that for all $r_1 \in R_1$ and all $x \in A^*:$ 
 \ $\varphi(r_1x) = \varphi(r_1) \, x$.
A {\it right ideal isomorphism} over $BA^*$ is a homomorphism, as above, 
such that the domain $R_1$ and the image $R_2$ are essential ideals, and 
such that $\varphi$ is bijective. 
Two right ideal homomorphisms $\psi: Q_1 \to Q_2$ and 
$\varphi: R_1 \to R_2$ are {\it essentially equal} iff $Q_1 =_{\sf ess} R_1$
and $\psi$ agrees with $\varphi$ on $Q_1 \cap R_1$; this implies that we
also have $Q_2 =_{\sf ess} R_2$.

Every right ideal homomorphism $\varphi$ over $BA^*$ has a {\it unique
maximal essentially equal extension} to a right ideal homomorphism of $BA^*$
(which is denoted ${\sf max}(\varphi)$). 
This can be proved in the same way as for right ideal homomorphisms over 
$A^*$ (see Prop.\ 1.2 in \cite{BiThomMon} and Prop.\ 2.1 in 
\cite{BiThomps}).
When $\varphi$ is an isomorphism, ${\sf max}(\varphi)$ is also an isomorphism.

To define $G_{k,i}$ we let the underlying set consist of all maximally 
extended right ideal isomorphisms between essential right ideals of $BA^*$.
The multiplication of $G_{k,i}$ is functional composition, followed by 
maximal extension (to a maximal right ideal isomorphism). This is similar 
to the definition of $G_{k,1}$ in \cite{BiThomps}; a similar definition
(with a different terminology) appears in \cite{Scott}.
We define the monoid $M_{k,i}$ by using maximally extended essentially 
equal right ideal homomorphisms between right ideals of $BA^*$. The 
multiplication is composition followed by maximal essentially equal 
extension. This is similar to the definition of $M_{k,1}$ in \cite{BiThomMon}.
Along similar lines one can define ${\it Inv}_{k,i}$, consisting of all
maximally extended essentially equal right ideal isomorphisms between (not 
necessarily essential) right ideals of $BA^*$. Compare with the definition 
of ${\it Inv}_{k,1}$ in \cite{BiThomMon}.
We do not need to assume that the alphabets $A$ and $B$ are disjoint.
We refer to Section 1 of \cite{BiThomMon} and Section 1 of \cite{BiRL} 
for terminology that is not defined here.

\medskip

\noindent
Here are some nice facts about $G_{k,i}$ (discovered by Higman \cite{Hig74},
see also \cite{ScottTour}):

\smallskip

\noindent $\bullet$ If $i \equiv j$ {\sf mod} $k-1$ then $G_{k,i}$ and 
$G_{k,j}$ are isomorphic (Coroll.\ 2, page 12 in \cite{Hig74}). 
So, in the notation ``$G_{k,i}$'' we can (and will) always assume that
$1 \leq i \leq k-1$. We will show that this holds for $M_{k,i}$ too. 

\smallskip

\noindent $\bullet$ By Theorem 6.4 in \cite{Hig74}: 
If $h \neq k$ then $G_{h,i}$ is not isomorphic to $G_{k,j}$ (for any 
$i, j$). 
Also, when ${\sf gcd}(k-1, i) \neq {\sf gcd}(k-1, j)$ 
then $G_{k,i}$ is not isomorphic to $G_{k,j}$. We will show that this holds 
for $M_{k,i}$ too. 
E.g., for all $1 \leq i \leq k-2$, \ $G_{k,i}$ is not isomorphic to 
$G_{k,k-1}$, and $M_{k,i}$ is not isomorphic to $M_{k,k-1}$.

\smallskip

\noindent $\bullet$ However (Theorem 7.3 in \cite{Hig74}), if $d$ divides 
$k$ then $G_{k,i}$ is isomorphic to $G_{k, di}$ (where $di$ is taken 
${\sf mod} \ k-1$). 
E.g., when $k$ is even, $G_{k,1}$ is isomorphic to $G_{k,2}$ and to 
$G_{k,k/2}$.
Hence by all these observation:  \ $G_{3,1} \not\simeq G_{3,2}$,  
 \ $G_{4,1} \simeq G_{4,2} \not\simeq G_{4,3}$.  

\smallskip

\noindent $\bullet$ Every group $G_{k,i}$ is finitely presented. When $k$ 
is even, $G_{k,i}$ is a simple group, and when $k$ is odd, $G_{k,i}$ 
contains a simple subgroup of index 2.

\medskip

We will show that the {\it maximal subgroups} of $M_{k,1}$ are isomorphic to
the Higman groups $G_{k,i}$ ($1 \leq i \leq k-1$). Thus, in $M_{k,1}$ we
``rediscover'' all the Higman groups $G_{k,i}$.

\smallskip

The monoid $M_{k,1}$ and the inverse monoid ${\it Inv}_{k,1}$ are finitely
generated \cite{BiThomMon}. 

Since $M_{k,1}$ acts partially on $A^*$, and in particular, $M_{2,1}$ acts 
partially on the set of all bit-strings $\{0,1\}^*$, we can view the 
elements of $M_{k,1}$ as boolean functions. 
In order to formalize this connection between $M_{k,1}$ and combinational 
boolean circuits we will also use an infinite generating set for 
$M_{k,1}$, of the
form $\Gamma \cup \tau$, where $\Gamma$ is any finite generating set of 
$M_{k,1}$, and $\tau$ consists of the letter position transpositions on 
strings. More precisely, \ $\tau = \{ \tau_{i,i+1} : i \geq 1 \}$, where
 \ $\tau_{i,i+1}(u \ x_i \ x_{i+1} \ v) \ = \ u \ x_{i+1} \ x_i \ v$, for
all $u \in A^{i-1}$, $v \in A^*$, and $x_i, x_{i+1} \in A$.

Then, for every combinational circuit $C$ there is a word $w$ over
$\Gamma \cup \tau$ such that: (1) \ the functions represented by $C$ and $w$
are the same, (2) \ $|w| \leq c \cdot |C|$ \ (for some constant $c$ which 
depends only on the choice of generators and gates). Here, $|C|$ is the
size of the circuit $C$ (i.e., the number of gates, plus the number of wire
crossings, plus the number of input or output ports), and $|w|$ is the 
length of the word $w$ over $\Gamma \cup \tau$; for this we define 
$|\tau_{i,i+1}| = i+1$ and $|\gamma| = 1$ for all $\gamma \in \Gamma$.  

Conversely, if a function \ $f: A^m \to A^n$ \ is represented by a
word $w$ over $\Gamma \cup \tau$ then $f$ has a combinational
circuit $C$ with \ $|C| \leq c \cdot |w|^2$ \ (for some constant $c$). 
See \cite{BiDistor}, Section 2.

We call a generating set of $M_{k,1}$ of the form $\Gamma \cup \tau$, as 
above, a {\it circuit-like generating set}.

\medskip

The {\it Green relations} $\leq_{\cal J}$, $\leq_{\cal L}$, $\leq_{\cal R}$, 
$\equiv_{\cal D}$, and $\leq_{\cal H}$ are classical concepts in the study 
of monoids (and semigroups), see e.g.\ \cite{CliffPres, Grillet}. 
By definition, for any $u,v \in M$ (where $M$ is a monoid) we have: 
$u \leq_{\cal J} v$ iff every ideal of $M$ 
containing $v$ also contains $u$; equivalently, $u \leq_{\cal J} v$ iff 
there exist $x,y \in M$ such that $u = xvy$.
Similarly, $u \leq_{\cal L} v$ iff any left ideal of $M$ containing $v$
also contains $u$; equivalently, there exists $x \in M$ such that $u = xv$;
the definition of $\leq_{\cal R}$ is similar.
By definition, $u \equiv_{\cal D} v$ iff there exists $s \in M$ such that
$u \equiv_{\cal R} s \equiv_{\cal L} v$; this is equivalent to the existence
of $t \in M$ such that $u \equiv_{\cal L} t \equiv_{\cal R} v$.
The $\cal H$-preorder is defined by $y \leq_{\cal H} x$ iff 
$y \leq_{\cal R} x$ and $y \leq_{\cal L} x$. 
For any pre-order $\leq_{\cal X}$ we define the corresponding equivalence 
relation $\equiv_{\cal X}$ by $y \equiv_{\cal X} x$ iff $y \leq_{\cal X} x$ 
and $x \leq_{\cal X} y$.

In \cite{BiThomMon} we gave characterizations of $\leq_{\cal J}$ and
$\equiv_{\cal D}$ in $M_{k,1}$. In \cite{BiRL} we characterized 
$\leq_{\cal L}$ and $\leq_{\cal R}$ in $M_{k,1}$, and we analyzed the 
computational complexity of deciding $\leq_{\cal L}$ or $\leq_{\cal R}$.

The main goal of this paper is to study the computational complexity of 
deciding $\equiv_{\cal D}$ and $\leq_{\cal J}$ in $M_{k,1}$. 
The problems of deciding whether $\psi \leq_{\cal J} \varphi$, or deciding 
whether $\psi \equiv_{\cal D} \varphi$, when $\psi$ and $\varphi$ are 
given by words over a finite generating set of $M_{k,1}$ 
(or of ${\it Inv}_{k,1}$), are in {\sf P}.
However, when the inputs $\psi$ and $\varphi$ are given by words over a 
circuit-like generating set, then deciding $\leq_{\cal J}$ for $M_{k,1}$ is
{\sf coDP}-complete, and deciding $\equiv_{\cal D}$ is 
$\oplus_{k-1} \! \bullet \! {\sf NP}$-complete. The complexity class {\sf DP} 
(called ``difference P''), introduced in \cite{PapadYannak}, has not been
used much in the literature; see Section 5 for details. The complexity class 
$\oplus_h \! \bullet \! {\sf NP}$ (for a given $h \geq 2$) is a counting 
complexity class; it fits into a pattern that has appeared in the 
literature; but this particular class has never been studied; see Section 
6 for details. 
There are related problems for circuits (see Sections 5 and 6) that are 
also complete for these unusual complexity classes.
In addition, we study the complexity of some search problems associated with
$\equiv_{\cal D}$ and $\leq_{\cal J}$ in $M_{k,1}$. 

We characterize the complexity of deciding the Green relations of 
${\it Inv}_{k,1}$. 
In particular, deciding whether $\psi \equiv_{\cal D} \varphi$ when $\psi$ 
and $\varphi$ are given by words over $\Gamma_I \cup \tau$ (where $\Gamma_I$
is a finite generating set of ${\it Inv}_{k,1}$), is 
$\oplus_{k-1} {\sf P}$-complete. The class $\oplus_h {\sf P}$ is a familiar
counting complexity class. For details, see Section 7.

%%%%%%%%%%%%%%%%%%%%%%%%%%%%%%%%%%%%%%%%%%%%%%%%%%%
%% Section 
%%%%%%%%%%%%%%%%%%%%%%%%%%%%%%%%%%%%%%%%%%%%%%%%%%%

\section{The maximal subgroups of $M_{k,1}$}

We saw in \cite{BiThomMon} (Prop.\ 2.2 and Theorem 2.5) that $M_{k,1}$ has 
only one non-zero $\cal J$-class, and that it has $k-1$ non-zero 
$\cal D$-classes. These $\cal D$-classes, denoted by $D_i$ for 
$i = 1, \ldots, k-1$, are given by

\medskip

$D_i \ = \ \{\varphi \in M_{k,1} \ : $
$ \ |{\sf imC}(\varphi)| \equiv i \ {\sf mod} \ k-1 \}$.

\medskip

\noindent  
It is well known and easy to see that every {\em subgroup of a semigroup} 
is an $\equiv_{\cal H}$-class, and that an $\equiv_{\cal H}$-class $H$ is 
a group iff $H$ contains an idempotent. 
The $\equiv_{\cal H}$-classes that contain an idempotent are the 
{\it maximal subgroups} of the semigroup, i.e., the subgroups that are 
not strictly contained in another subgroup.
It is well known and not hard to prove that all maximal subgroups of a 
same $\cal D$-class are isomorphic (see e.g.\ \cite{Grillet}  Prop.\ 2.1
and the remark that follows it, or \cite{Lallement} Coroll.\ 2.7).

We saw (\cite{BiThomps} Prop.\ 2.1) that $G_{k,1}$ is the group of units 
of $M_{k,1}$ (i.e., the group of invertible elements). This implies that 
$M_{h,1}$ is not isomorphic to $M_{k,1}$ when $h \neq k$ (since we know
from \cite{Hig74} that $G_{h,1}$ is not isomorphic to $G_{k,1}$ when 
$h \neq k$).
The fact that $M_{k,1}$ has $k-1$ non-zero $\cal D$-classes also implies
$M_{h,1} \not\simeq M_{k,1}$ when $h \neq k$.
 
The next theorem shows a very nice correspondence between the $k-1$ 
non-zero $\cal D$-classes and the $k-1$ groups $G_{k,i}$ 
($1 \leq i \leq k-1$) that Higman introduced in \cite{Hig74}. It is 
surprising (at first) that all the $G_{k,i}$ show up automatically in the 
structure of $M_{k,1}$.  

\begin{thm} \label{max_subgroups_Mk1} 
 \  For every $i$ ($1 \leq i \leq k-1$) we have: The maximal subgroups of 
the $\cal D$-class $D_i$ of $M_{k,1}$ are isomorphic to the Thompson-Higman 
group $G_{k,i}$. 
\end{thm}
{\bf Proof.} In the $\cal D$-class $D_i$ we consider the idempotent 
$\eta_i  = {\sf id}_{\{a_1, \ldots, a_i\}}$, i.e., the partial identity map 
that is defined on those (and only those) words that start with a 
letter in $\{a_1, \ldots, a_i\}$. Since $|{\sf imC}(\eta_i)| = i$ 
we have indeed $\eta_i \in D_i$. Consider the set

\smallskip

$G_{\eta_i} \ = \ $
$\{ \varphi \in M_{k,1} : {\sf Dom}(\varphi)$ and ${\sf Im}(\varphi)$ 
are essential right subideals of \ $\{a_1, \ldots, a_i\} \, A^* \}$.

\smallskip

\noindent The set $G_{\eta_i}$ is a subgroup of $M_{k,1}$, with identity 
element $\eta_i$. Moreover, this group is isomorphic to $G_{k,i}$; an 
isomorphism is obtained by replacing each $b_j w \in B \, A^*$ 
($1 \leq j \leq i$) by $a_j w \in A \, A^*$. 
Clearly, the subgroup $G_{\eta_i}$ is contained in the $\cal H$-class
of $\eta_i$. 

Conversely, suppose $\varphi \equiv_{\cal H} \eta_i$.  Then 
$\varphi$ is injective with domain essentially equal to 
$\{a_1, \ldots, a_i\} \, A^*$ (since $\varphi \equiv_{\cal L} \eta_i$, and
by the characterization of $\leq_{\cal L}$ in Section 3.4 of \cite{BiRL}).
And the image of $\varphi$ is essentially equal to 
$\{a_1, \ldots, a_i\} \, A^*$ (since $\varphi \equiv_{\cal R} \eta_i$, and
by the characterization of $\leq_{\cal R}$ in Section 2 of \cite{BiRL}). 
It follows that $\varphi \in G_{\eta_i}$, by the definition of 
$G_{\eta_i}$.   So $G_{\eta_i}$ is the entire $\equiv_{\cal H}$-class of
$\eta_i$, hence it is a maximal subgroup, in $D_i$. Since all the 
maximal subgroups in the same $\equiv_{\cal D}$-class $D_i$ are isomorphic,
every maximal subgroup of $M_{k,1}$ is isomorphic to some $G_{\eta_i}$ 
(which is itself isomorphic to $G_{k,i}$). 
 \ \ \ $\Box$

%%%%%%%%%%%%%%%%%%%%%%%%%%%%%%%%%%%%%%%%%%%%%%%%%%%%%%%%%%%%%%
% Section
%%%%%%%%%%%%%%%%%%%%%%%%%%%%%%%%%%%%%%%%%%%%%%%%%%%%%%%%%%%%%%

%%%%%
\section{ The Thompson-Higman monoids $M_{k,i}$ }

In the Introduction we defined $M_{k,i}$ by using two alphabets, 
$A = \{a_1, \ldots, a_k\}$, and $B_i = \{b_1, \ldots, b_i\}$. 
It follows from this definition that when $1 \leq j \leq i$, 
$M_{k,j}$ is a submonoid of $M_{k,i}$ (not just up to isomorphism, but 
also as a subset).

The identity element of $M_{k,i}$ can be described by the table
${\sf id}_{B_i} = \{(b,b) : b \in B_i\}$, and will be denoted by {\bf 1} 
(if $k$ and $i$ are clear from the context).

\begin{pro} \label{MkiMkjwheniequivj} \  
If $s \equiv t$ mod $k-1$ then $M_{k,s} \simeq M_{k,t}$.
\end{pro}
{\bf Proof.} It suffices to prove that for all $n \geq k$ we have
$M_{k,n} \simeq M_{k,n-(k-1)}$.  We embed $B_nA^*$ into $B_{n-k+1}A^*$ by 
the map 

\smallskip
 
$E :  \ \ \left\{ \begin{array}{ll}   
b_i \ \longmapsto \ b_i \ \ & {\rm for} \ i = 1, \ldots, n-k; \\  
b_{i+n-k} \ \longmapsto \ b_{n-k+1} a_i \ \ & {\rm for} \ i = 1, \ldots, k.
\end{array} \right. $

\smallskip

\noindent The image of this embedding is the essential right ideal \, 
$B_{n-k}A^* \, \cup \, b_{n-k+1}A A^*$, which is an essential right 
sub-ideal of $B_{n-k+1} A^*$.
The embedding $B_nA* \hookrightarrow B_{n-k+1}A^*$ determines an
embedding $M_{k,n} \hookrightarrow M_{k,n-(k-1)}$ that we will also call
$E$. The embedding is surjective since it is the identity on the submonoid 
$M_{k,n-(k-1)}$ of $M_{k,n}$; hence the embedding is also a retract.  

The embedding is a homomorphism: Consider any $\psi, \varphi \in M_{k,1}$. 
After essential restrictions, if needed, we can assume that 
$\varphi, \psi \in M_{k,n}$ have tables of the form 
$\{(u_i,v_i) : i \in I\}$, respectively $\{(v_j,w_j): j \in J\}$, such that 
the set $\{v_i : i \in I\} \cup \{v_j: \in J\}$ is a prefix code. 
So the product $\psi \varphi(.)$ is represented by the composition of these 
tables (without need to extend or restrict), i.e., $\psi \varphi(.)$ has a 
table $\{(u_i,w_i): i \in I \cap J\}$.
Then  $E(\psi) \cdot E(\varphi)(.)$ has a table 
$\{(E(u_i),E(w_i)): i \in I \cap J\}$. On the other hand, by applying $E$ to
the table $\{(u_i,w_i): i \in I \cap J\}$ for $\psi \varphi(.)$, we see that 
$\{(E(u_i),E(w_i)): i \in I \cap J\}$ is also a table for
$E(\psi \cdot \varphi)(.)$; hence, $E(\psi \cdot \varphi)(.)$
$ = $ $E(\psi) \cdot E(\varphi)(.)$.
 \ \ \ $\Box$

\medskip

  From now on, when we write $M_{k,i}$ we will always assume that 
$1 \leq i \leq k-1$.

\medskip

By definition, the group of units of a monoid $M$ is the set of invertible 
elements of $M$; equivalently, the group of units is the maximal subgroup 
of $M$ whose identity is the identity of the monoid. 
\begin{pro} \label{groupunitsMki} \ 
The group of units of $M_{k,i}$ and of ${\it Inv}_{k,i}$ is $G_{k,i}$.
\end{pro}
{\bf Proof.} The proof is very similar to the proof of Prop.\ 2.1 in 
\cite{BiThomMon} (which shows that $G_{k,1}$ is the group of units of 
$M_{k,1}$).
 \ \ \ $\Box$

\begin{cor} \label{ifGnonoisoThenMnoniso} \  
If $G_{k,i} \not\simeq G_{h,j}$ then $M_{k,i} \not\simeq M_{h,j}$.
\end{cor}
{\bf Proof.} If two monoids have non-isomorphic groups of units then they 
are non-isomorphic.
 \ \ \ $\Box$

\begin{thm} \label{MkiGreenrelations} {\bf (Green relations of 
$M_{k,s}$).} \ For all $k \geq 2$ and $s$ (with $1 \leq s \leq k-1$) we
have:   \\  
{\bf (1)} \ The monoids $M_{k,s}$ and ${\it Inv}_{k,s}$ are 
0-$\cal J$-simple, i.e., they have only one non-zero $\cal J$-class. \\ 
{\bf (2)} \ $M_{k,s}$ and ${\it Inv}_{k,s}$ are congruence-simple. \\ 
{\bf (3)} \ For all $\psi, \varphi \in M_{k,s}$ (or ${\it Inv}_{k,s}$):
 \ \ $\varphi \equiv_{\cal D} \psi$ \ \ iff 
 \ \ $|{\sf imC}(\varphi)| \equiv |{\sf imC}(\psi)| \ {\sf mod} \ k-1$.

 \ \ Hence $M_{k,s}$ and ${\it Inv}_{k,s}$ have $k-1$ non-zero 
  $\cal D$-classes. \\     
{\bf (4)} The $\leq_{\cal R}$ and $\leq_{\cal L}$ preorders for $M_{k,s}$ 
have the same characterizations as for $M_{k,1}$ 

 \ \ (Theorems 2.1 and 3.32 in \cite{BiRL}).
\end{thm}
{\bf Proof.} {\bf (1)} \ The proof of Prop.\ 2.2 in \cite{BiThomMon} can
easily adapted to $M_{k,s}$ and ${\it Inv}_{k,s}$.

When $\varphi \in M_{k,s}$ (or $\in {\it Inv}_{k,s}$) is not the empty map, 
there exist $b_m x_0, b_n y_0 \in BA^*$ such that 
$b_n y_0 = \varphi(b_m x_0)$. 
Let $P = \{p_1, \ldots, p_s\} \subset A^*$ be a prefix code (not necessarily
maximal) with $|P| = s$.
Then we have $\varphi(b_m x_0 p_i) = b_n y_0 p_i$ for $i = 1, \ldots, s$.
Let us define $\alpha, \beta \in {\it Inv}_{k,s}$ by the tables 
$\alpha = \{(b_i, b_m x_0 p_i) : i = 1, \ldots, s\}$ and 
$\beta = \{(b_n y_0 p_i, b_i): i = 1, \ldots, s\}$.
Then $\beta \varphi \alpha(.) = \{(b_i, b_i) : i = 1, \ldots, s\} = {\bf 1}$. 
So every non-zero element of $M_{k,s}$ (or ${\it Inv}_{k,s}$) is in the same 
$\cal J$-class as the identity element.

\smallskip

\noindent
{\bf (2)} \ The proof of congruence-simplicity is exactly the same as the 
proof of Theorem 2.3 in in \cite{BiThomMon}.
 
\smallskip

\noindent
{\bf (3)} \ The proof of Theorem 2.5 in \cite{BiThomMon} works for $M_{k,s}$
and ${\it Inv}_{k,s}$ too. Proposition 2.4 in \cite{BiThomMon} remains 
unchanged, and Lemma 2.6 becomes: 

{\it For all finite alphabets $A$, $B$, and 
every integer $i \geq 0$ there exists a maximal prefix code in $BA^*$ of 
cardinality $|B| + (|A| - 1)i$. And every finite maximal prefix code in 
$BA^*$ has cardinality $|B| + (|A| - 1)i$ for some integer $i \geq 0$.
}

\smallskip

\noindent Lemma 2.7 remains unchanged. 

\smallskip

\noindent
Let $A = \{a_1, \ldots, a_k\}$ and $B = \{b_1, \ldots, b_s\}$, be the finite 
alphabets used in the definition of $M_{k,s}$.
The statement of Lemma 2.8 becomes: \\
(1) For any $m \geq k + s - 1$, let $i$ be the residue of $m-(s-1)$ modulo 
$k-1$ in the range $2 \leq i \leq k$, and let us write 
$m = s-1 + i + (k-1)j$, for some $j \geq 0$. Then there exists a prefix code 
$Q'_{i,j}$ of cardinality $|Q'_{i,j}| = m$, such that 
${\sf id}_{Q'_{i,j}}$ is an essential restriction of 
${\sf id}_{\{b_1, \ldots, b_{s-1}, b_s a_1, \ldots, b_s a_i\}}$.
Hence ${\sf id}_{Q'_{i,j}} = $
$ {\sf id}_{\{b_1, \ldots, b_{s-1}, b_s a_1, \ldots, b_s a_i\}}$ as
elements of ${\it Inv}_{k,s}$. \\   
(2) In $M_{k,s}$ and in ${\it Inv}_{k,s}$ we have 
 \ $ {\sf id}_{\{b_1, \ldots, b_{s-1}, b_s a_1\}} $ $\equiv_{\cal D} $
$ {\sf id}_{\{b_1, \ldots, b_{s-1}, b_s a_1, \ldots, b_s a_k\}} $ $=$
${\bf 1}$. 

\smallskip

\noindent In the proof of Lemma 2.8(1), $Q_{i,j}$ is replaced by
 \ $Q'_{i,j} = \{b_1, \ldots, b_{s-1}\} \ \cup \ b_1 \, Q_{i,j}$. 

\smallskip

\noindent Lemmas 2.9 and 2.10 are unchanged.

\smallskip

\noindent
In the final proof of Theorem 2.5 we replace the end of the first paragraph
by the following:
In particular, when $|Q_1| \equiv s-1+i$ {\sf mod} $k-1$ (with
$1 \leq i \leq k$), then $\varphi_1 \equiv_{\cal D} $
${\sf id}_{\{b_1, \ldots, b_{s-1}, b_s a_1, \ldots, b_s a_i\}}$. 

\smallskip

\noindent
{\bf (4)} \ The proofs of Theorems 2.1 and 3.32 in \cite{BiRL} are
straightforwardly generalized to $M_{k,s}$. 
We will see later in Proposition 7.1 that (4) also holds for 
${\it Inv}_{k,1}$; and for ${\it Inv}_{k,s}$ the proof is similar.
 \ \ \ $\Box$

\begin{pro} \label{maxsubGrMki} \ 
The maximal subgroups of $M_{k,i}$ are isomorphic to $G_{k,j}$ for 
$j = 1, \ldots, k-1$, with $G_{k,j}$ being isomorphic to the
maximal subgroup of the $\cal D$-class \ $D_j = $
$\{\varphi \in M_{k,i} : |{\sf imC}(\varphi)| \equiv j \ {\sf mod} \ k-1\}$.
The same is true for ${\it Inv}_{k,i}$.
\end{pro}
{\bf Proof.} This is similar to the proof of Theorem \ref{max_subgroups_Mk1} 
above. In the $\cal D$-class $D_j$ we can pick, for example, the idempotent
${\sf id}_{\{b_1, \ldots, b_j\}}$ if $1 \leq j < s$, and we can pick the
idempotent
${\sf id}_{\{b_1, \ldots, b_{s-1}, b_s a_1, \ldots, b_s a_{j - s +1}\}}$
if $s \leq j \leq k-1$.
 \ \ \ $\Box$

\medskip

Since the Green relations ${\cal J, D, R, L}$ of $M_{k,i}$ are quite
similar to those of $M_{k,1}$, we will focus on $M_{k,1}$ from now on.

%%%%%%%%%%%%%%%%%%%%%%%%%%%%%%%%%%%%%%%%%%%%%%%%%%%%%%%%%%%%%%
% Section
%%%%%%%%%%%%%%%%%%%%%%%%%%%%%%%%%%%%%%%%%%%%%%%%%%%%%%%%%%%%%%

%%%%%
\section{Complexity of \, $\leq_{\cal J}$ \, and \, $\equiv_{\cal D}$ 
 \, over a finite generating set }

We are interested in the difficulty of checking on input 
$\psi, \varphi \in M_{k,1}$, whether $\psi \leq_{\cal J} \varphi$, or 
$\psi \equiv_{\cal J} \varphi$, or $\psi \equiv_{\cal D} \varphi$.
In \cite{BiRL} we addressed the question whether 
$\psi \leq_{\cal R} \varphi$, or $\psi \leq_{\cal L} \varphi$.
We assume at first that $\psi, \varphi \in M_{k,1}$ are given either by
tables, or by words over a chosen finite generating set $\Gamma$ of 
$M_{k,1}$. Recall that $M_{k,1}$ is finitely generated (Theorem 3.4 in 
\cite{BiThomps}). For computational complexity it does not matter much 
which finite generating set of $M_{k,1}$ is used; finite changes in the 
generating set only lead to linear changes in the complexity.

Let {\bf 0} denote the zero element of $M_{k,1}$ (represented by the empty
map), and let {\bf 1} denote the identity element of $M_{k,1}$ (represented
by the identity map on $A^*$).

Checking whether $\psi \leq_{\cal J} \varphi$ is not difficult. Since
$M_{k,1}$ is 0-$\cal J$-simple (Prop.\ 2.2  in \cite{BiThomMon}) we have:
 \ $\psi \leq_{\cal J} \varphi$ iff $\varphi \neq {\bf 0}$ or
$\psi = {\bf 0}$. Consider any element $\psi \in M_{k,1}$, given by a 
table or by a word over a chosen finite generating set $\Gamma$ of $M_{k,1}$.
In order to check whether $\psi$ is equal to {\bf 0}, we calculate
${\sf imC}(\psi)$, as an explicit list of words. If $\psi$ is given by
a table, ${\sf imC}(\psi)$ can be directly read from the table.
If $\psi$ is given by a word over a finite generating set of $M_{k,1}$ we
use Corollary 4.11 of \cite{BiThomMon} to find the list of elements of 
${\sf imC}(\psi)$ in polynomial time. To check whether $\psi = {\bf 0}$ 
we now check whether ${\sf imC}(\psi) = \varnothing$.

The relation $\psi \equiv_{\cal D} \varphi$ can be checked in deterministic
polynomial time, by using the characterization of $\equiv_{\cal D}$ in
Theorem 2.5 in \cite{BiThomMon} 
(which says that \ $\psi \equiv_{\cal D} \varphi$ \ iff 
 \ $|{\sf imC}(\psi)| \equiv |{\sf imC}(\varphi)|$ {\sf mod} $k-1$).
We can compute ${\sf imC}(\psi)$ and
${\sf imC}(\varphi)$ as explicit lists of words, either from the table or
by Corollary 4.11 of \cite{BiThomMon}, in polynomial time.

This proves:

\begin{pro} \label{J_D_order_probl_fg} \
The $\leq_{\cal J}$ decision problem and the $\equiv_{\cal D}$ decision
problem of $M_{k,1}$ are decidable in deterministic polynomial time,
if inputs are given by tables or by words over a finite generating set.
 \ \ \ $\Box$
\end{pro}

In connection with the $\leq_{\cal J}$-relation we consider the
{\bf multiplier search problem} for $M_{k,1}$ over a finite generating set
$\Gamma$.  This problem is specified as follows: \\
{\sf Input:} \ $\varphi, \psi \in M_{k,1}$, given by words over $\Gamma$. \\
{\sf Premise:} \ $\psi \leq_{\cal J} \varphi$. \\   
{\sf Search}: \ Find some $\alpha, \beta \in M_{k,1}$, given by words over 
$\Gamma$, such that \ $\psi = \beta \varphi \alpha(.)$.

Note that since the decision problem for $\leq_{\cal J}$ over a finite set 
of generators is in {\sf P}, the premise is easily checked, so this is 
problem could be reformulated without a premise.  

\begin{pro} \label{Jorder_search-fg} \
The $\leq_{\cal J}$-relation multiplier search problem for $M_{k,1}$ is
solvable in deterministic polynomial time, if inputs and output are given by
words over a finite generating set.
\end{pro}
{\bf Proof.} If $\psi = {\bf 0}$ we pick $\alpha = \beta = {\bf 0}$.
Let us assume now that $\psi \neq {\bf 0} \neq \varphi$.
We can choose the multipliers $\alpha, \beta \in M_{k,1}$ as 
follows (as we did already in the proof of 0-$\cal J$-simplicity, i.e., 
Proposition 2.2 in \cite{BiThomMon}).   

First, from $\varphi$ (given by a word over $\Gamma$) we want to find 
some $x_0, y_0 \in A^*$ such that $y_0 = \varphi(x_0)$; we want to do 
this in deterministic polynomial time (as a function of 
$|\varphi|_{\Gamma}$). 
By Corollary 4.11 in \cite{BiThomMon} we find an explicit list
of ${\sf imC}(\varphi)$ in polynomial time. In this list we pick any 
element $y_0 \in {\sf imC}(\varphi)$. From $y_0$ and the generator sequence
for $\varphi$ we can then find an element $x_0 \in \varphi^{-1}(y_0)$
as follows.  By Corollary 4.15 in \cite{BiThomMon} we can, in deterministic
polynomial time, build a deterministic partial finite automaton that 
accepts the set $\varphi^{-1}(y_0)$. By a search in this finite automaton 
we can (in deterministic polynomial time) find a word $x_0$ that is accepted 
by the automaton.

Now let \ $\alpha = \{(\varepsilon, x_0)\}$ \ and 
 \ $\beta' = \{(y_0, \varepsilon)\}$. Since $\alpha$ and $\beta'$ have tables 
with one entry of polynomial length, we can (in polynomial time) find words 
over $\Gamma$ that represent $\alpha$, respectively $\beta'$; for this we 
use Lemma 5.3 of \cite{BiRL} (which, in polynomial time, finds a word over
$\Gamma$ from a table). 

Now we have \ $\beta' \varphi \alpha(.) = {\bf 1}$. Hence
 \ $\psi \beta' \varphi \alpha(.) = \psi$. Clearly, a word over $\Gamma$ 
for $\beta = \psi \beta'$ can be found in deterministic polynomial time, 
since we can find a word for $\beta'$ in deterministic polynomial time. 
 \ \ \ $\Box$

\medskip

In connection with the $\equiv_{\cal D}$-relation we consider the 
$\cal D$-{\bf pivot search problem} for $M_{k,1}$ over a finite generating 
set $\Gamma$.  This problem is specified as follows: \\ 
{\sf Input:} \ $\varphi, \psi \in M_{k,1}$, given by words over $\Gamma$. \\  
{\sf Premise:} \ $\varphi \equiv_{\cal D} \psi$. \\  
{\sf Search:} \ Find an element $\chi \in M_{k,1}$, given by a word over 
$\Gamma$, such that $\psi \equiv_{\cal R} \chi \equiv_{\cal L} \varphi$.  

Note that since the decision problem for $\equiv_{\cal D}$ over a finite 
set of generators is in {\sf P}, the premise is easily checked, so this  
problem can easily be transformed to an ordinary search problem, without 
premise.

\begin{pro} \label{D_pivot_search-fg} \  
The $\cal D$-pivot search problem for the $\equiv_{\cal D}$-relation of 
$M_{k,1}$ is solvable in deterministic polynomial time, if inputs and output 
are given by words over a finite generating set. 
\end{pro}
{\bf Proof.} As in the problem statement, let $\varphi, \psi \in M_{k,1}$ 
with $\varphi \equiv_{\cal D} \psi$; so,
$|{\sf imC}(\varphi)| \equiv |{\sf imC}(\psi)|$ {\sf mod} $k-1$. 
By Corollary 4.11 in \cite{BiThomMon}, ${\sf imC}(\varphi)$ and  
${\sf imC}(\psi)$ can be found in deterministic polynomial time (and hence
they have polynomial size). In a polynomial number of steps, we can 
essentially restrict $\varphi$ and $\psi$ to $\varphi'$, respectively 
$\psi'$ such that $|{\sf imC}(\varphi')| = |{\sf imC}(\psi')|$.
We can obtain the restricted map $\varphi'$ by taking \ $\varphi'$ 
$ = {\sf id}_{{\sf imC}(\varphi')} \circ \varphi(.)$. Since 
$|{\sf imC}(\varphi')|$ is polynomially bounded in terms of 
$|\varphi|_{\Gamma}$, the map ${\sf id}_{{\sf imC}(\varphi')}$ has a
polynomially bounded table, and hence a word over $\Gamma$ can be found for 
${\sf id}_{{\sf imC}(\varphi')}$ in polynomial time (by Lemma 5.2 in
\cite{BiRL}).
Thus we obtain a word over $\Gamma$ for $\varphi'$ in polynomial time, 
and similarly for $\psi'$.
Let $\alpha$ be any element of $M_{k,1}$ that maps ${\sf imC}(\psi')$ 
bijectively onto ${\sf imC}(\varphi')$. Since ${\sf imC}(\psi')$ and
${\sf imC}(\varphi')$ can be explicitly listed in polynomial time, we can
find a table (and hence a word over $\Gamma$, by Lemma 5.2 in \cite{BiRL}) 
for $\alpha$ in polynomial time.  Then we have:

\smallskip

$\psi' \ \equiv_{\cal L} \ \alpha \psi' \ \equiv_{\cal R} \ \varphi'$ .

\smallskip 

\noindent The latter $\equiv_{\cal R}$ holds because $\alpha \psi'$ 
is a map from ${\sf domC}(\psi')$ onto ${\sf imC}(\varphi')$, hence 
$\alpha \psi'$ and $\varphi'$ have the same image code (which implies 
$\equiv_{\cal R}$ by Theorem 2.1 of \cite{BiRL}).
Thus, $\alpha \psi'$ is a $\cal D$-pivot. Since $\alpha$ and $\psi'$ can be 
found in deterministic polynomial time, we can find a word for this 
$\cal D$-pivot in deterministic polynomial time.
 \ \ \ $\Box$

%%%%%%%%%%%%%%%%%%%%%%%%

\section{The complexity of \, $\leq_{\cal J}$ \, over the generating set 
            $\Gamma \cup \tau$ }

We consider the $\leq_{\cal J}$ decision problem and the $\leq_{\cal J}$
multiplier search problem of $M_{k,1}$ over the {\it circuit-like 
generating set} $\Gamma \cup \tau$, where $\Gamma$ is any chosen finite 
generating set of $M_{k,1}$, and $\tau = \{\tau_{i,i+1} : i \geq 1\}$. 
As we saw near the end of the Introduction, this generating set makes the
elements of $M_{k,1}$ similar to combinational circuits: Circuit-size 
becomes polynomially equivalent to the word-length 
\cite{BiCoNP,BiDistor,BiThomMon,BiFact}. 
The word problem and the Green relations of $M_{k,1}$ over $\Gamma$ are
in {\sf P}. But over $\Gamma \cup \tau$ the word problem of $M_{k,1}$ is 
{\sf coNP}-complete \cite{BiThomMon}, the $\leq_{\cal R}$ decision problem 
is $\Pi_2^{\sf P}$-complete, and the $\leq_{\cal L}$ decision problem 
is {\sf coNP}-complete \cite{BiRL}. 

For complexity and word-length, finite changes in the generating set do 
not matter much; they only lead to linear changes in the complexity or
the word-length. So, for a circuit-like generating set
$\Gamma \cup \tau$ we can choose $\Gamma$ arbitrarily, provided that 
$\Gamma$ is finite and $\Gamma \cup \tau$ generates $M_{k,1}$.

%%%

\subsection{The $\leq_{\cal J}$ decision problem over $\Gamma \cup \tau$ }

Because of the 0-$\cal J$-simplicity of the $\cal J$-order we have to 
consider the following special word problem in $M_{k,1}$ over
$\Gamma \cup \tau$. 

\smallskip

\noindent {\sf Input:} \ $\varphi \in M_{k,1}$, given by a word over
the generating set $\Gamma \cup \tau$,

\smallskip

\noindent {\sf Question} ({\bf 0 word problem}):
 \ Is \ $\varphi = {\bf 0}$ \ as an element of $M_{k,1}$ ?

\medskip

\noindent Recall that the word problem in $M_{k,1}$ over 
$\Gamma \cup \tau$ is {\sf coNP}-complete (Theorem 4.12 in \cite{BiThomMon}). 
In \cite{BiRL} (Prop.\ 6.2) we proved the following:

\begin{pro} \label{0_word_probl_infgen} \
The {\bf 0} word problem of $M_{k,1}$ over $\Gamma \cup \tau$ is
{\sf coNP}-complete.
\end{pro}
{\bf Proof.} We reduce the tautology problem for boolean formulas to the
{\bf 0} word problem. Let $B$ be any boolean formula, with corresponding
boolean function $\{0,1\}^m \to \{0,1\}$. We identify $\{0,1\}$ with
$\{a_1,a_2\} \subseteq \{a_1, \ldots, a_k\} = A$.
The function $B$ can be viewed as an element $\beta \in M_{k,1}$,
represented by a word over $\Gamma \cup \tau$. The length of that word is
linearly bounded by the size of the formula $B$ (by Prop.\ 2.4 in
\cite{BiDistor}). In $M_{k,1}$ we consider the element ${\sf id}_{0 A^*}$
(i.e., the identity function restricted to $0A^*$), and we assume that
some fixed representation of ${\sf id}_{0 A^*}$ by a word over $\Gamma$
has been chosen. We have:

\smallskip

 \ \ \ ${\sf id}_{0 A^*} \circ \beta(.) = {\bf 0}$
 \ \ iff \ \ ${\sf Im}(\beta) \subseteq  1 \, A^*$.

\smallskip

\noindent The latter holds iff $B$ is a tautology. Thus we reduced the
tautology problem for $B$ to the special word problem
 \ ${\sf id}_{0A^*} \, \beta = {\bf 0}$. Note that ${\sf id}_{0A^*}$ is
fixed, and independent of $B$.

It follows that the {\bf 0} word problem of $M_{k,1}$ over
$\Gamma \cup \tau$ is {\sf coNP}-hard for all $k \geq 2$. Moreover,
since the word problem of $M_{k,1}$ over $\Gamma \cup \tau$ is in 
{\sf coNP} (by Prop.\ 4.12 in \cite{BiThomMon}), it follows that the 
{\bf 0} word problem is {\sf coNP}-complete.
 \ \ \ $\Box$.

\bigskip

We can now characterize the complexity of the decision problem
of the $\cal J$-order of $M_{k,1}$ over $\Gamma \cup \tau$.
We will need the following complexity classes:

\medskip

${\sf DP} \ = \ {\sf NP} \wedge {\sf coNP} \ = \ $
$\{ L_1 \cap L_2 : L_1 \in {\sf NP} \ {\sf and} \ L_2 \in {\sf coNP}\}$
$ \ = \ $
$\{ N_1 - N_2 \ : \ N_1, N_2 \in {\sf NP}\}$,

\medskip

${\sf coDP} \ = \ {\sf NP} \vee {\sf coNP} \ = \ $
$\{ L_1 \cup L_2 : L_1 \in {\sf NP} \ {\sf and} \ L_2 \in {\sf coNP}\}$.

\medskip

\noindent In other words, {\sf DP} consists of the set-differences between
pairs of sets in {\sf NP}.  The class {\sf DP} was introduced
in \cite{PapadYannak}, where several problems were proved to be
{\sf DP}-complete (see also pp.\ 92-95 in \cite{Handb}). 
In particular, the following problem, called {\sf Sat-and-unsat} was 
given as an example of a {\sf DP}-complete problem:
The input consists of two boolean formulas $B_1$ and $B_2$, and the
question is whether $B_1$ is satisfiable and $B_2$ is unsatisfiable.
It follows immediately that the following problem is also
{\sf DP}-complete; the input is as before, and the question is
whether $B_1$ is not a tautology {\em and} $B_2$ is a tautology.
Hence, the following problem, which we call
{\sf Nontaut-or-taut}, is {\sf coDP}-complete: \\
{\sf Input:} \ Two boolean formulas $B_1$ and $B_2$. \\
{\sf Question:} \ Is $B_1$ is not a tautology {\em or} is $B_2$ a 
tautology? \ (I.e., 
 \ $(\forall x_1) B_1(x_1) \stackrel{?}{\Rightarrow} (\forall x_2) B_2(x_2)$)

%Equivalently, is \ $[(\forall \vec{x}_1) B_1(\vec{x}_1)$
%$\Rightarrow$ $(\forall \vec{x}_2) B_2(\vec{x}_2)]$ \ true?
%Here, $\vec{x}_1$ and $\vec{x}_2$ denote the lists of boolean variables of
%$B_1$ respectively $B_2$.

\smallskip

The class {\sf coDP} is closed under union and under polynomial-time
disjunctive reduction, whereas {\sf DP} is closed under intersection and
under polynomial-time conjunctive reduction.
The classes {\sf DP} and {\sf coDP} constitute the second level of the
{\em boolean hierarchy} {\sf BH}; for more information on {\sf DP} and
{\sf BH}, see e.g.\ the survey \cite{CGHHSWW}.

\begin{thm} \label{J_order_infgen} \ In $M_{k,1}$ over the generating set
$\Gamma \cup \tau$ we have:  

\smallskip

\noindent {\bf (1)} \ The $\equiv_{\cal J} {\bf 0}$ decision problem is 
  {\sf coNP}-complete. 

\smallskip

\noindent {\bf (2)} \ The $\equiv_{\cal J} {\bf 1}$ decision problem is 
  {\sf NP}-complete. 

\smallskip

\noindent {\bf (3)} \ The $\equiv_{\cal J}$ and $\leq_{\cal J}$ decision 
problems for $M_{k,1}$ over $\Gamma \cup \tau$ are ${\sf coDP}$-complete 
(for the $\equiv_{\cal J}$ decision problem, the ${\sf coDP}$-completeness 
is with respect to polynomial-time disjunctive reductions).
\end{thm}
{\bf Proof. (1)} \ In any semigroup, $\equiv_{\cal J} {\bf 0}$ is 
equivalent to $= {\bf 0}$. We saw that the {\bf 0} word problem is 
{\sf coNP}-complete (Prop.\ \ref{0_word_probl_infgen} above).  

\smallskip

\noindent
{\bf (2)} \ By 0-$\cal J$-simplicity of $M_{k,1}$, \ 
$\varphi \equiv_{\cal J} {\bf 1}$ iff $\varphi \neq {\bf 0}$. So, the
$\equiv_{\cal J} {\bf 1}$ decision problem is equivalent to the negation
of the {\bf 0} word problem, hence it is {\sf NP}-complete. 

\smallskip

\noindent
{\bf (3)} \ The $\equiv_{\cal J}$- and $\leq_{\cal J}$-decision problems 
are in ${\sf coNP} \vee {\sf NP}$ because (by 0-$\cal J$-simplicity of 
$M_{k,1}$), $\psi \leq_{\cal J} \varphi$ is equivalent to $\psi = {\bf 0}$
or $\varphi \neq {\bf 0}$ as elements of $M_{k,1}$.
The question whether $\psi = {\bf 0}$ is in {\sf coNP}, and the question
whether $\varphi \neq {\bf 0}$  is in {\sf NP}.

Let us prove ${\sf coDP}$-hardness of the $\leq_{\cal J}$ decision
problem.  For boolean formulas $B_1$ and $B_2$ we have: \, $B_1$ is not a 
tautology or $B_2$ is a tautology \ iff
 \ ${\sf id}_{0\{0,1\}^*} \circ \beta_1 = {\bf 0}$ \ or
 \ ${\sf id}_{0\{0,1\}^*} \circ \beta_2 \neq {\bf 0}$, which is iff
${\sf id}_{0\{0,1\}^*} \circ \beta_1$ $\leq_{\cal J} $
${\sf id}_{0\{0,1\}^*} \circ \beta_2$. This reduces the 
{\sf Nontaut-or-taut} problem to the $\leq_{\cal J}$ decision problem.

The $\equiv_{\cal J}$ decision problems is ${\sf coDP}$-hard because 
the $\leq_{\cal J}$ decision problem reduces to it by a polynomial-time
disjunctive reduction: $\psi \leq_{\cal J} \varphi$ \ iff 
 \ $\psi \equiv_{\cal J} {\bf 0}$ or $\psi \equiv_{\cal J} \varphi$.
The class {\sf coDP} is closed under union and under polynomial-time
disjunctive reduction.  
 \ \ \ $\Box$.

%%%%

\subsection{The $\leq_{\cal J}$ multiplier search problem}

The multiplier search problem for $M_{k,1}$ over $\Gamma \cup \tau$ is 
specified as follows: \\
{\sf Input:} \ $\varphi, \psi \in M_{k,1}$, given by words over 
$\Gamma \cup \tau$. \\
{\sf Premise:} \ $\psi \leq_{\cal J} \varphi$. \\
{\sf Search}: \ Find some $\alpha, \beta \in M_{k,1}$, expressed as words 
over $\Gamma \cup \tau$, such that \ $\psi = \beta \varphi \alpha(.)$.

By 0-$\cal J$-simplicity of $M_{k,1}$ the multiplier search problem is
trivial when $\psi$ or $\varphi$ are {\bf 0}. 
When $\psi$ and $\varphi$ are not {\bf 0}, both will be 
$\equiv_{\cal J} {\bf 1}$.

Therefore we consider the {\bf special multiplier search problem} for 
$\equiv_{\cal J} {\bf 1}$ in $M_{k,1}$ over $\Gamma \cup \tau$, specified 
as follows: \\ 
{\sf Input:} \ $\varphi \in M_{k,1}$, given by a word over 
 $\Gamma \cup \tau$. \\  
{\sf Premise:} \ ${\bf 1} \equiv_{\cal J} \varphi$. \\  
{\sf Search:} \ Find one $\alpha$ and one $\beta \in M_{k,1}$, described 
by words over $\Gamma \cup \tau$, such that 
$\beta \, \varphi \, \alpha = {\bf 1}$.
 
When we have multipliers $\alpha$ and  $\beta$ such that 
${\bf 1} = \beta \varphi \alpha$ then we can take the pair
$\psi \beta$, $\alpha$ to obtain multipliers for 
$\psi \leq_{\cal J} \varphi$.
 
We saw (Prop.\ \ref{Jorder_search-fg}) that the problem is solvable in 
deterministic polynomial time when $M_{k,1}$ is taken over any finite 
generating set $\Gamma$. 
Note that over $\Gamma \cup \tau$, the premise (namely that
${\bf 1} \equiv_{\cal J} \varphi$) is non-trivial, being {\sf NP}-complete.

\smallskip

See the Appendix for general information on search problems, the classes 
{\sf NPsearch} and {\sf xNPsearch}, search reductions, and completeness.

\smallskip

Before we deal with the multiplier search problem for $\leq_{\cal J}$
we will consider the {\bf domain element search problem} of $M_{k,1}$ 
over $\Gamma \cup \tau$. The problem is specified as follows. \\ 
{\sf Input:} \ $\varphi \in M_{k,1}$, given by a word over 
  $\Gamma \cup \tau$. \\ 
{\sf Premise:} \ $\varphi \neq {\bf 0}$. \\ 
{\sf Search:} \ Find an element $x_0 \in {\sf Dom}(\varphi)$. \\ 
The corresponding relation is 
 \ $\{(\varphi, x_0) \in (\Gamma \cup \tau)^* \times A^* \ : \ $
$\varphi(x_0) \neq \varnothing \}$.

A similar problem is the {\bf inverse image search problem} of $M_{k,1}$
over $\Gamma \cup \tau$, specified as follows. \\  
{\sf Input:} \ $y_0 \in A^*$, and $\varphi \in M_{k,1}$, given by a word over
  $\Gamma \cup \tau$. \\
{\sf Premise:} \ $y_0 \in {\sf Im}(\varphi)$. \\
{\sf Search:} \ Find an element $x_0 \in \varphi^{-1}(y_0)$.

\begin{pro} \label{Domainwitness_search-infg} \ 
The domain element search problem and the inverse image search problem 
of $M_{k,1}$ over $\Gamma \cup \tau$ are {\sf xNPsearch}-complete.
\end{pro}
{\bf Proof.} \ The longest words in ${\sf domC}(\varphi)$ have length 
$\leq c \cdot |\varphi|_{\Gamma \cup \tau}$, for some constant
$c$ (by Theorem 4.5 in \cite{BiThomMon}); the constant $c$ is the length 
of the longest word in the tables of the elements of $\Gamma$.
Hence there exists $x_0 \in {\sf domC}(\varphi)$ with polynomial length,
and in fact, all elements of ${\sf domC}(\varphi)$ have polynomial length. 
So, without loss of existence of solutions, we can consider the 
polynomially balanced sub-problem 

\smallskip 

 $\{(\varphi, x_0) \ : \ \varphi(x_0) \neq \varnothing \ \ {\rm and} \ $
 $x_0 \in {\sf domC}(\varphi) \}$.
 
\smallskip

\noindent By Prop.\ 5.5 in \cite{BiRL}, we can verify in deterministic 
polynomial time whether $x_0 \in {\sf domC}(\varphi)$. Hence, this
sub-problem is in {\sf NPsearch}.

In order to reduce the {\sc SatSearch} problem to the domain element
search problem of $M_{k,1}$ over $\Gamma \cup \tau$, we can view a boolean 
circuit $B$ as an element of $M_{k,1}$, given by a word over 
$\Gamma \cup \tau$. Then $x_0$ is an element of the domain of
${\sf id}_{1\{0,1\}^*} \circ B(.)$ iff $x_0$ satisfies $B$. 

Essentially the same proof works for the inverse image search problem. 
 \ \ \ $\Box$

\begin{pro} \label{Jorder_search-infg} \   
The special multiplier search problem for $\equiv_{\cal J} {\bf 1}$ in 
$M_{k,1}$ over $\Gamma \cup \tau$ is {\sf xNPsearch}-complete. 
In particular, for any $\varphi \equiv_{\cal J} {\bf 1}$ there exist 
multipliers of polynomial word-length over $\Gamma \cup \tau$. 
\end{pro}
{\bf Proof.} To show that the problem is in {\sf xNPsearch} we follow the
proof of Prop.\ \ref{Jorder_search-fg} above and of Proposition 2.2 in
\cite{BiThomMon}. When $\varphi \neq {\bf 0}$ there exists 
$x_0 \in {\sf domC}(\varphi)$; let  $y_0 = \varphi(x_0)$.
Since the longest words in ${\sf domC}(\varphi) \cup {\sf imC}(\varphi)$
have length $\leq c \cdot |\varphi|_{\Gamma \cup \tau}$, for some constant
$c$ (by Theorem 4.5 in \cite{BiThomMon}), $x_0$ and $y_0$ have polynomial
length. Then 
 \ $\beta \varphi \alpha = \{(\varepsilon,\varepsilon)\} = {\bf 1}$, where 
 \ $\alpha = \{(\varepsilon, x_0)\}$ \ and
 \ $\beta = \{(y_0, \varepsilon)\}$. 
Thus, we take the sub-problem defined by the following relation:

\smallskip

$\{ \big(\varphi, (\{(y_0, \varepsilon)\}, \ \{(\varepsilon, x_0)\})\big)$
  $ \ : \ \varphi \in (\Gamma \cup \tau)^*, \ x_0 \in {\sf domC}(\varphi),$
  $ \ y_0 = \varphi(x_0) \}$.

\smallskip

\noindent We saw that 
$|x_0|, |y_0| \leq c \cdot |\varphi|_{\Gamma \cup \tau}$; hence this 
relation is polynomially balanced. 
The verification problem for this relation is in {\sf P}: Indeed, we 
can check in deterministic polynomial time whether 
$x_0 \in {\sf domC}(\varphi)$ (by Prop.\ 5.5 in \cite{BiRL}). We can
compute $\varphi(x_0)$ in deterministic polynomial time (by the proof of
Theorem 4.12 in \cite{BiThomMon}), and compare $\varphi(x_0)$ with $y_0$.

Since $\alpha$ and $\beta$ have tables with one entry of polynomial size, 
we can (in polynomial time) find words over $\Gamma$ that represent 
$\alpha$, respectively $\beta$; for this we use Lemma 5.3 of \cite{BiRL} 
(which, in polynomial time, finds a word over $\Gamma$ from a table). 
So, $x_0$ and $y_0$ yield multipliers (expressed as strings over
$\Gamma \cup \tau$) for $\varphi \equiv_{\cal J} {\bf 1}$.

\smallskip

To show {\sf NPsearch}-completeness we reduce the problem {\sc SatSearch}
to the $\equiv_{\cal J} {\bf 1}$ multiplier search problem over 
$\Gamma \cup \tau$. (See the Appendix for the definition of search 
reductions.)
We construct an input-output reduction $(\rho_{\rm in}, \, \rho_{\rm sol})$ 
as follows.
The function $\rho_{\rm in}$ maps any boolean formula $B(x_1, \ldots, x_m)$
to ${\sf id}_{1\{0,1\}^*} \circ B(.) \in M_{k,1}$;
here, ${\sf id}_{1\{0,1\}^*}$ is the partial identity    
with domain and image $1\{0,1\}^*$. From the boolean formula for $B$ we 
easily construct a word over $\Gamma \cup \tau$ for $B$; moreover, 
we can choose a fixed word over $\Gamma$ to represent 
${\sf id}_{1\{0,1\}^*}$.
For all $t \in \{0,1\}^m$, \ $w \in \{0,1\}^*$:

\smallskip

 \ \ \ ${\sf id}_{1\{0,1\}^*} \circ B(tw) \ = \ $
$\left\{ \begin{array} {ll}
  1w   &  \mbox{if \ $B(t) = 1$,} \\   
  \varnothing &  \mbox{otherwise.}
          \end{array} \right.  $

\smallskip

\noindent We have
 \ ${\bf 1} \equiv_{\cal J} {\sf id}_{1\{0,1\}^*} \circ B$
 \ iff there are multipliers $\beta, \alpha$ such that 
 \ ${\bf 1} = \beta \circ {\sf id}_{1\{0,1\}^*} \circ B \circ \alpha$. 
By what we saw, these multipliers can be chosen as follows:
 \ $\beta = \{(y_0, \varepsilon)\}$, 
 \ $\alpha = \{(\varepsilon, x_0)\}$, with 
 \ $x_0 \in {\sf domC}({\sf id}_{1\{0,1\}^*} \circ B)$, and
 \ $y_0 = {\sf id}_{1\{0,1\}^*} \circ B(x_0)$.
Moreover, ${\sf domC}({\sf id}_{1\{0,1\}^*} \circ B) = $
$\{ t \in \{0,1\}^m : B(t) = 1\}$, hence for 
$x_0 \in {\sf domC}({\sf id}_{1\{0,1\}^*} \circ B)$ we have $y_0 = 1$. 

Therefore, the multiplier $\alpha = \{(\varepsilon, x_0)\}$ determines 
a solution of {\sc SatSearch}, by reading $x_0$ in the table of $\alpha$.
So, we simply define the map $\rho_{\rm sol}$ by
 \ $\rho_{\rm sol}(\beta, \alpha) = \alpha(\varepsilon)$ ($= x_0$).

Finally, to obtain a verification reduction we consider the map 
 \ $\rho_{\rm ver} : (B,t) \longmapsto $
$({\sf id}_{1\{0,1\}^*} \circ B, \ \beta, \ \alpha)$, where 
 \ $\beta = \{(1, \varepsilon)\}$, and \ $\alpha = \{(t, \varepsilon)\}$. 
Then $\rho_{\rm ver}$ reduces the verification problem 
``$B(t) \stackrel{?}{=} 1$'' to the verification problem  
 \ ``$\beta \circ {\sf id}_{1\{0,1\}^*} \circ B \circ \alpha$
$\stackrel{?}{=} {\bf 1}$''.
Indeed,
 \ $\beta \circ {\sf id}_{1\{0,1\}^*} \circ B \circ \alpha = $
$\{(\varepsilon,\varepsilon)\}  = {\bf 1}$ 
 \ iff 
 \ $\beta \circ {\sf id}_{1\{0,1\}^*} \circ B \circ \alpha(\varepsilon)$
$ \ = \ \beta \circ {\sf id}_{1\{0,1\}^*} \circ B(t) \ = \ \varepsilon$,
which holds when $B(t) = 1$, and does not hold when $B(t) = 0$. 
 \ \ \ $\Box$

%%%

\subsection{The multiplier search problems for $\leq_{\cal R}$ and 
$\leq_{\cal L}$}  

By definition, a {\it left- (right-) inverse} of an element $x$ in a monoid 
$M$ is an element $t \in M$ such that $tx = {\bf 1}$ (respectively 
$xt = {\bf 1}$). 
A {\it left multiplier for} $\psi \leq_{\cal L} \varphi$ in $M$ is any 
$\beta \in M$ such that $\psi = \beta \varphi$. A {\it right multiplier for}
$\psi \leq_{\cal R} \varphi$ in $M$ is any $\alpha \in M$ such that 
$\psi = \varphi \alpha$.

\begin{pro} \label{multSearchLR} \hspace{-.12in} {\bf .} 
For $M_{k,1}$ over $\Gamma \cup \tau$ we have: \\
{\bf (1)} \ The left multiplier search problem for $\leq_{\cal L}$, and the 
left-inverse search problem are {\em not} in 
{\sf xNPsearch} (unless {\sf NP} = {\sf coNP}).  \\
{\bf (2)} \ The right multiplier search problem for $\leq_{\cal R}$, 
and the right-inverse search problem are {\em not} in 
{\sf xNPsearch} (unless {\sf NP} = {\sf coNP}).    
\end{pro}
{\bf Proof.} \ (1) \ If the problems were in {\sf xNPsearch} we could guess
a polynomial-size multiplier, and for some such guess the verification 
problem would be in {\sf P} (by the definition of {\sf xNPsearch}).
Hence, the $\leq_{\cal L}$ and $\equiv_{\cal L} {\bf 1}$ decision problems
would be in {\sf NP}.
However we saw in \cite{BiRL} (Section 6.2) that these two problems are
{\sf coNP}-complete. Hence, we would have {\sf NP} = {\sf coNP}, i.e., the
polynomial hierarchy would collapse to level 1.

(2) \ If the search problems were in {\sf xNPsearch} then (by the same 
reasoning as for the $\cal L$-order)
the $\leq_{\cal R}$ and $\equiv_{\cal R} {\bf 1}$ decision problems
would be in {\sf NP}. However, we saw in  \cite{BiRL} (Section 2.2) that 
these two problems are $\Pi_2^{\sf P}$-complete. Hence we would have 
$\Pi_2^{\sf P} = {\sf NP}$, hence {\sf coNP} = {\sf NP} (since 
$\Pi_2^{\sf P}$ contains {\sf coNP}).  \ \ \ $\Box$

\medskip

Note that we also saw in \cite{BiRL} (Section 5.3) that for $M_{k,1}$ over
$\Gamma \cup \tau$ we have:  Unless the polynomial hierarchy collapses to 
level 2, the $\leq_{\cal R}$-multipliers and the right-inverses
do not have polynomially bounded word-length.

%%%%%%%%%%%%%%%%%%%%%%%%%%%%%%%%%%%%%

\section{The complexity of \, $\equiv_{\cal D}$ \, over the generating set
            $\Gamma \cup \tau$ }

Recall the characterization of the $\cal D$-relation of $M_{k,1}$: 
There are $k-1$ non-{\bf 0} $\cal D$-classes, $D_1, \ldots, D_{k-1}$;  
for any $\varphi \in M_{k,1}$ we have \ $\varphi \in  D_i$ \ iff 
 \ $|{\sf imC}(\varphi)| \equiv  i  \ \ {\sf mod} \ \ k-1$.
Since $M_{k,1}$ has only finitely many $\cal D$-classes, the 
$\equiv_{\cal D}$-decision problem is equivalent to the membership 
problems of these $k$ $\cal D$-classes.

The $\cal D$-class $\{ {\bf 0} \}$ is special. 
Membership in the $\cal D$-class $\{ {\bf 0} \}$ is the same thing as  
the {\bf 0} word problem, which is {\sf coNP}-complete.  

In order to characterize the complexity of the membership problem of a
non-zero $\cal D$-class $D_i$ we need a somewhat exotic complexity class. 

%%%

\subsection{New counting complexity classes}

Recall Valiant's counting complexity class $\# {\sf P}$ (pronounced 
``number P''), consisting of all functions $f_R: A^* \to \{0,1\}^*$ of 
the form

\smallskip

 \ \ \ \ \  $f_R(x) \ = \ $ binary representation of the number 
 \ $|\{ y \in A^* : (x,y) \in R \}|$ ;

\smallskip

\noindent here $R$ ranges over all predicates 
$R \subseteq A^* \times A^*$ such that the membership problem 
``$(x,y) \stackrel{?}{\in} R$'' is in {\sf P}
(deterministic polynomial time), and such that $R$ is polynomially
balanced. A predicate $R$ is called {\it polynomially balanced} iff 
there exists a polynomial $p$ such that for all $(x,y) \in R$: 
$|y| \leq p(|x|)$;  see e.g.\ p.\ 181 in \cite{Papadim}, and note that 
the definition of ``balanced'' is not symmetric in $x$ and $y$.

This can be generalized: In the above definition we replace {\sf P} by 
any complexity class $\cal C$; then we obtain a counting class 
$\# \bullet {\cal C}$, corresponding to polynomially balanced predicates 
whose membership problem is in $\cal C$. For the history of these complexity
classes, and in particular, for the reason why there is a dot in the notation,
see  \cite{TodaDiss, HemaspVollmer, DurandHerKol}. 
The classes $\# \bullet {\sf NP}$ and $\# \bullet {\sf coNP}$ have been 
studied and, in particular, is was proved that 
 \ $\# \bullet {\sf NP} \ = \ \# \bullet {\sf coNP}$ 
 \ iff \ ${\sf NP} = {\sf coNP}$ \ \cite{KoSchoTo}. 

Another important counting class is $\oplus {\sf P}$, introduced in
\cite{PapadZachos} and \cite{GP86}. More generally, $\oplus_{h,i} {\sf P}$ 
consists of all sets $L_R$ of the form

\smallskip

 \ \ \ \ \  $L_R \ = \ \big\{x \in A^* \ : \ $
$|\{ y \in A^* : (x,y) \in R \}| \equiv  i  \ {\sf mod} \ h \big\}$ ;

\smallskip

\noindent here $R$ ranges over all predicates $R \subseteq A^* \times A^*$
such that the membership problem ``$(x,y) \stackrel{?}{\in} R$'' is in 
{\sf P}, and such that $R$ is polynomially balanced. 
And $h,i$ are integers with $h \geq 2$.  In that notation, $\oplus {\sf P}$
is $\oplus_{2,1} {\sf P}$.  It was proved that if 
 \ $\oplus {\sf P} \subseteq \Sigma_{\ell}^{\sf P}$ \ then the polynomial
hierarchy collapses to \ $\Sigma_{\ell+1}^{\sf P} \cap \Pi_{\ell+1}^{\sf P}$ 
 \ (due to Toda \cite{Toda91}; see also \cite{DuKo} pp.\ 334-340).
The notation {\sc mod}$_h {\sf P}$ was used in 
\cite{CaiHem90, BeiGillHert, BeiGill} for \ ${\sf co} \oplus_{h,0} {\sf P}$
 \ ($ = \bigcup_{i=1}^{h-1} \oplus_{h,i} {\sf P}$).
See pp.\ 297-298 of \cite{HemOgiCompan} for some properties of 
$\oplus {\sf P}$ and {\sc mod}$_h {\sf P}$.

The class $\oplus_{h,i} {\sf P}$ can be generalized to
$\oplus_{h,i} \bullet {\cal C}$ for any class $\cal C$ of formal 
languages, and in particular to $\oplus_{h,i} \! \bullet \! {\sf NP}$. 
The class $\oplus_{2,0} \bullet {\cal C}$
was mentioned in \cite{HemaspVollmer}. 
For a predicate $R \subseteq A^* \times A^*$ and any $x_1,x_2 \in A^*$ we 
use the notation 

\smallskip

 $(x_1)R = \{x_2 \in A^* : (x_1,x_2) \in R\}$, \ and

\smallskip
 
 \ $R(x_2) = \{x_1 \in A^* : (x_1,x_2) \in R\}$. 

\medskip

For a predicate $R \subseteq A^* \times A^*$ we say that $R \in {\cal C}$ 
\ iff \ the language $\{ xby \in A^*bA^*: (x,y) \in R\}$ belongs to 
${\cal C}$, for some letter $b \not\in A$.

\begin{defn} \label{oplus_h_i_dotC} \     
Let $h \geq 2$ and $i \geq 0$. 
A set $L \subseteq A^*$ belongs to $\oplus_{h,i} \bullet {\cal C}$ iff 
there is a polynomially balanced predicate $R \in {\cal C}$ such that for 
all $x \in A^*$: 
 
\smallskip

  $x \in L$ \ \ iff \ \ $|(x)R| \equiv  i  \ {\sf mod} \ h$. 

\smallskip

\noindent In that case we say that $L$ can be defined (in 
$\oplus_{h,i} \bullet {\cal C}$) by the predicate $R$. 
\end{defn}
Note that (except when $h = 2$) this definition is unsymmetric for 
$x \in L$ versus $x \not\in L$; so when $h > 2$, 
 \ $\oplus_{h,i} \bullet {\cal C}$ and 
${\sf co} \oplus_{h,i} \bullet {\cal C}$ seem to be different (but this
remains an open question).
Obviously, if $i \equiv j$ {\sf mod} $h$ then 
$\oplus_{h,i} \bullet {\cal C} = \oplus_{h,j} \bullet {\cal C}$.

\begin{lem} \label{h_i_Add} \   
Let $L$ be defined in $\oplus_{h,i} \! \bullet \! {\sf NP}$ by a predicate 
$R$. Then $L$ can also be defined in $\oplus_{h,i+1} \! \bullet \! {\sf NP}$ 
by a predicate $R'$, such that for all $x \in A^*$: \ $|(x)R'| = |(x)R| + 1$. 
\end{lem}
{\bf Proof.} Let $L \in \oplus_{h,i} \! \bullet \! {\sf NP}$ be defined by 
a polynomially balanced predicate $R \in {\sf NP}$. Let us denote 
$\{(x,x) : x \in A^* \}$ by $\Delta$, and let us assume for the moment that
$\Delta \cap R = \varnothing$. Then we have: 

\smallskip

$L \ = \ \big\{x \in A^* \ : \ $
$|\{ y \in A^* : (x,y) \in R \cup \Delta \}| \equiv i+1 $
$ \ {\sf mod} \ h \big\}$.

\smallskip

\noindent Indeed, for any $x \in A^*$: \ $(x)(R \cup \Delta) = $
$(x)R \cup \{x\}$, \ and \ $x \not\in (x)R$ \ since 
$\Delta \cap R = \varnothing$.
So \ $|(x)(R \cup \Delta)| \equiv i+1 \ {\sf mod} \ h$. 
The predicate $R \cup \Delta$ is in {\sf NP}, and it is polynomially
balanced. Thus, $L$ is also defined  (in 
$\oplus_{h,i+1} \! \bullet \! {\sf NP}$) by the predicate $R \cup \Delta$.

If $R$ does not satisfy $\Delta \cap R = \varnothing$, we consider 
 \ $R' = \{ (x,xya_1) \in A^* \times A^* : (x,y) \in R\}$, which is
polynomially balanced and in {\sf NP}, and satisfies
$\Delta \cap R' = \varnothing$; here, $a_1$ is one of the letters of $A$. 
Moreover, $R'$ defines $L$ as a element of 
$\oplus_{h,i} \! \bullet \! {\sf NP}$ since \  
$|\{ y \in A^* : (x,y) \in R \}| = $
$|\{ y \in A^* : (x, xya_1) \in R' \}| = $
$|\{ z \in A^* : (x, z) \in R' \}|$. 
Now we can replace $R$ by $R'$ and carry out the previous reasoning, which 
assumed that $\Delta \cap R' = \varnothing$.
 \ \ \ $\Box$

\medskip 

\noindent By applying Lemma \ref{h_i_Add} at most $h-1$ times we obtain: 

\begin{cor} \label{h_i_AllSame} \
For all $i, j:$ \ \ $\oplus_{h,i} \! \bullet \! {\sf NP} \ = \ $
$\oplus_{h,j} \! \bullet \! {\sf NP}$. I.e., for any fixed $h \geq 2$, the
classes  $\oplus_{h,i} \! \bullet \! {\sf NP}$ are the same for all $i$.
\end{cor}
Therefore we will use the notation $\oplus_h \! \bullet \! {\sf NP}$ for each 
$\oplus_{h,i} \! \bullet \! {\sf NP}$. 
 
\begin{lem} \label{NPcoNPsubPlusNP}  
 \ Both {\sf NP} and {\sf coNP} are subsets of 
$\oplus_h \! \bullet \! {\sf NP}$.
Moreover, every $L_0$ in ${\sf NP} \cup {\sf coNP}$ can be defined (in
$\oplus_{h,1} \! \bullet \! {\sf NP}$) by a predicate $R$ such that: 

\smallskip

$x \in L_0$ \ \ iff 
 \ \ $|\{ y : (x,y) \in R\}| \equiv 1 \ {\sf mod} \ h$, \ \ \ and 

\smallskip

$x \not\in L_0$ \ \ iff 
 \ \ $|\{ y : (x,y) \in R\}| \equiv 0\ {\sf mod} \ h$.
\end{lem}
{\bf Proof.} Let $L \in {\sf NP}$, let $x \in A^*$, and let 
$\{u_1, \ldots, u_{h-1}\} \subset A^*$ be a fixed set of $h-1$ different 
non-empty words.  Then we have

\smallskip

$\{ y \in A^* \ : \ y = x \ {\sf and} \ x \in L\} \ = \ $
$\left\{  \begin{array}{ll} 
          \{x\}       &  \mbox{if \ $x \in L$ ,} \\   
          \varnothing &  \mbox{if \ $x \not\in L$ .}  
          \end{array}  \right. $

\smallskip
 
\noindent Similarly we have

\smallskip

$\big\{ y \in A^* \ : \ y = x \ {\sf or} $
$(y \in \{xu_1, \ldots, xu_{h-1}\} \ {\sf and} \ x \in L)\big\}$ $ \ = \ $
$\left\{  \begin{array}{ll}
 \{x\}       &  \mbox{if \ $x \in \overline{L}$ ,} \\
 \{x, xu_1, \ldots, xu_{h-1}\} &  \mbox{if \ $x \in L$ . }
\end{array}  \right. $

\smallskip
  
\noindent Hence, 

\smallskip

 \ \ \  $L \ = \ \big\{ x \in A^* \ : \ $
            $|\{ y \in A^* : y = x \ {\sf and} \ x \in L\}| \equiv 1 $
            {\sf mod} $h \big\}$,

\medskip

 \ \ \  $\overline{L} \ = \ \big\{ x \in A^* \ : \ $
            $|\{ y \in A^* : y = x \ {\sf or} $
              $(y \in \{xu_1, \ldots, xu_{h-1}\} \ {\sf and} \ x \in L)\}|$
            $\equiv 1 \ {\sf mod} \ h \big\}$.

\smallskip

\noindent  The predicate $R$ defined by  \ $(x,y) \in R$ \ iff 
 \ $[y = x \ {\sf and} \ x \in L]$, \ belongs to {\sf NP}, and is 
polynomially balanced.
Similarly, the predicate $R'$ defined by  \ $(x,y) \in R'$ \ iff \   
$[y = x \ {\sf or} $
$(y \in \{xu_1, \ldots, xu_{h-1}\} \ {\sf and} \ x \in L)]$ \ belongs to 
{\sf NP}, and is polynomially balanced.
So $L$ and $\overline{L}$ belong to $\oplus_{h,1} \! \bullet \! {\sf NP}$. 

One sees immediately form the definition of the predicates $R$ and $R'$ that
they have the following property: \ If $x \not\in L$
then $|\{ y \in A^* : (x,y) \in R\}| \equiv 0 \ {\sf mod} \ h$; 
if $x \not\in \overline{L}$ then 
$|\{ y \in A^* : (x,y) \in R'\}| = h \equiv 0 \ {\sf mod} \ h$. 
 \ \ \ $\Box$

\medskip

\noindent Lemma \ref{NPcoNPsubPlusNP} inspires the following definition.

\begin{defn} \label{defOplus_h_i_j} \    
For any integer $h \geq 2$ and two disjoint sets 
$S_1, S_2 \subset \{0, 1, \ldots, h-1\}$ we define the class 
$\oplus_{h,S_1,S_2} \! \bullet \! {\sf NP}$ as follows: 
 \ $L \subseteq A^*$ belongs to $\oplus_{h,S_1,S_2} \! \bullet \! {\sf NP}$
 \ iff \ there exists a polynomially balanced predicate 
$R \subseteq A^* \times A^*$ in {\sf NP} such that for all $x \in A^*$, 

\smallskip

$x \in L$ \ \ iff \ \ $|(x)R| \in S_1$ {\sf mod} $h$ , and 

\smallskip

$x \not\in L$  \ \ iff \ \ $|(x)R| \in S_2$ {\sf mod} $h$.

\smallskip 

\noindent We say then that $L$ can be defined in 
$\oplus_{h,S_1,S_2} \! \bullet \! {\sf NP}$ by the predicate $R$. 
In this notation the class $\oplus_{h,i} \! \bullet \! {\sf NP}$ is
$\oplus_{h,\{i\},\{j:j\neq i\}} \! \bullet \! {\sf NP}$.

When $S_1 = \{i\}$, $S_2 = \{j\}$ with $i \neq j$ we write
$\oplus_{h,i,j} \! \bullet \! {\sf NP}$.
\end{defn}
The second sentence of Lemma \ref{NPcoNPsubPlusNP} says that {\sf NP} and
{\sf coNP} are subclasses of $\oplus_{h,1,0} \! \bullet \! {\sf NP}$.

Clearly, \ {\sf co} $\oplus_{h,S_1,S_2} \! \bullet \! {\sf NP}$ $=$ 
$\oplus_{h,S_2,S_1} \! \bullet \! {\sf NP}$ \ (always assuming
$S_1 \cap S_2 = \varnothing$). 

By Lemma \ref{h_i_Add}, 
 \ $\oplus_{h,i,j} \! \bullet \! {\sf NP}$
$ = \oplus_{h,i+1,j+1} \! \bullet \! {\sf NP}$, 
and \ $\oplus_{h,S_1,S_2} \! \bullet \! {\sf NP} = $
$\oplus_{h,S_1+1,S_2+1} \! \bullet \! {\sf NP}$. By definition,
$S_1 + 1 = \{ i+1 : i \in S_1\}$, and similarly for $S_2 + 1$; all numbers
are taken modulo $h$. 

\begin{lem} \label{multiplying_i_j}  \   
Suppose $m$ is prime with $h$, and suppose that $L$ can be defined in 
$\oplus_{h,i,j} \! \bullet \! {\sf NP}$ by a predicate $R$. Then $L$ can also 
be defined in $\oplus_{h, m i, m j} \! \bullet \! {\sf NP}$ by a predicate $R'$ 
such that for all $x \in A^*$: 

\smallskip

$|(x)R'| = m \cdot |(x)R|$.

\smallskip

\noindent Hence,  
 \  $\oplus_{h,i,j} \! \bullet \! {\sf NP} = $
$ \oplus_{h, m i, m j} \! \bullet \! {\sf NP}$. (Here the numbers $m i$ and
$m j$ are taken modulo $h$.)
\end{lem}
{\bf Proof.}  By assumption, $x \in L$ iff $|(x)R| \equiv i$ {\sf mod} 
$h$, and $x \not\in L$ iff $|(x)R| \equiv j$ {\sf mod} $h$.
We choose a prefix code $\{u_1, \ldots, u_m\} \subset A^*$
of size $m$, and for $s = 1, \ldots, m$ we let 
 \ $R_s = \{ (x, u_s y) : (x,y) \in R\}$. Then each $R_s$ is in {\sf NP} 
and polynomially balanced. Moreover, $L$ is also defined in 
$\oplus_{h,i,j} \! \bullet \! {\sf NP}$ by $R_s$, since $|(x)R| = |(x)R_s|$.
Let \ $R' = R_1 \cup \ldots \cup R_m$. Then  $R'$ is also in {\sf NP} and
it is polynomially balanced. Since $\{u_1, \ldots, u_m\}$ is a prefix code 
we have $R_s \cap R_t = \varnothing$ when $s \neq t$. It follows that  
we have \ $|(x)R'| = m \cdot |(x)R|$ for all $x \in A^*$; and we have
$x \in L$ iff $|(x)R'| \equiv mi$ {\sf mod} $h$, and $x \not\in L$ iff 
$|(x)R'| \equiv mj$ {\sf mod} $h$.

Finally, when $m$ is prime with $h$ then $i \not\equiv j$ {\sf mod} $h$ 
implies $mi \not\equiv mj$  {\sf mod} $h$.  Hence 
$\oplus_{h,i,j} \! \bullet \! {\sf NP} = $
$ \oplus_{h, m i, m j} \! \bullet \! {\sf NP}$.
 \ \ \ $\Box$

\begin{cor} \label{allOplus_h_i_j} \   
For all $i, j$ such that $i - j$ is prime with $h$ we have:
 \  $\oplus_{h,i,j} \! \bullet \! {\sf NP} = $
$\oplus_{h,1,0} \! \bullet \! {\sf NP}$.
\end{cor}
{\bf Proof.} By $h-j$ applications of Lemma \ref{h_i_Add} we obtain 
 \ $\oplus_{h,i,j} \! \bullet \! {\sf NP} = $
$\oplus_{h,\ell,0} \! \bullet \! {\sf NP}$, where $\ell = i-j$ {\sf mod} $h$. 
Since $\ell$ is prime with $h$, $\ell$ has a mutiplicative inverse
$\ell^{-1}$ modulo $h$, hence by Lemma \ref{multiplying_i_j} (with 
$m = \ell^{-1}$) we obtain \ $\oplus_{h.\ell,0} \! \bullet \! {\sf NP} = $ 
$\oplus_{h,1,0} \! \bullet \! {\sf NP}$.
 \ \ \ $\Box$

\begin{cor} \label{Oplus_h_1_0_closComple} \ 
The class $\oplus_{h,1,0} \! \bullet \! {\sf NP}$ is closed under complement,
and contains {\sf NP} and {\sf coNP}.
\end{cor}
{\bf Proof.} We noted already that 
{\sf co} $\oplus_{h,1,0} \! \bullet \! {\sf NP}$ 
$ = \oplus_{h,0,1} \! \bullet \! {\sf NP}$. By Corollary \ref{allOplus_h_i_j},
$\oplus_{h,0,1} \! \bullet \! {\sf NP} = $
$\oplus_{h,1,0} \! \bullet \! {\sf NP}$. 
By Lemma \ref{NPcoNPsubPlusNP} and by Definition \ref{defOplus_h_i_j}, 
$\oplus_{h,1,0} \! \bullet \! {\sf NP}$ contains {\sf NP} and {\sf coNP}.
  \ \ \ $\Box$

\begin{lem} \label{Oplus_h_1_0_disjPred} \   
Suppose $L_1, L_2 \in \oplus_{h,1,0} \! \bullet \! {\sf NP}$ can be defined 
(in $\oplus_{h,1,0} \! \bullet \! {\sf NP}$) by predicates $R_1$, 
respectively $R_2$.  Then $L_1$ and $L_2$ can also  be defined (in 
$\oplus_{h,1,0} \! \bullet \! {\sf NP}$) by predicates $R_1'$, respectively 
$R_2'$ such that $R_1' \cap R_2' = \varnothing$.
\end{lem}
{\bf Proof.} We choose a prefix code $\{u_1, u_2\} \subset A^*$, and for
$i = 1,2$ we let \ $R_i' = \{ (x, u_i y) : (x,y) \in R\}$.
Then $R_1' \cap R_2' = \varnothing$ since $\{u_1, u_2\}$ is a prefix code.
Also, $R_i'$ is in {\sf NP} and is polynomially balanced.
And $|(x)R_i| = |(x)R_i'|$.  
  \ \ \ $\Box$

\begin{cor} \label{Oplus_h_1_0_closUnion} \     
For $h \geq 3$, if \ $L_1, L_2 \in \oplus_{h,1,0} \! \bullet \! {\sf NP}$ 
 \ then \ $L_1 \cap L_2 \in \oplus_h \! \bullet \! {\sf NP}$.

Hence for $h \geq 3$, \ ${\sf DP} \subseteq \oplus_h \! \bullet \! {\sf NP}$.
\end{cor}
{\bf Proof.} Let $R_1, R_2$ be predicates (in {\sf NP}) that describe 
$L_1$, respectively $L_2$ in $\oplus_{h,1,0} \! \bullet \! {\sf NP}$.
By Lemma \ref{Oplus_h_1_0_disjPred} we can assume that 
$R_1 \cap R_2 = \varnothing$.
So for all $x \in A^*:$ $(x)R_1 \cap (x)R_2 = \varnothing$,
hence \ $|(x)R_1 \cup (x)R_2| = |(x)R_1| + |(x)R_2|$. 

Since $L_1, L_2 \in \oplus_{h,1,0} \! \bullet \! {\sf NP}$ we have 
$|(x)R_1| \equiv$ 0 or 1 {\sf mod} $h$, according as $x \in L_1$ or
$x \not\in L_1$; and similarly for $L_2$.
Therefore, $x \in L_1 \cup L_2$ iff $|(x)(R_1 \ {\sf or} \ R_2)| \equiv$
1 or 2 {\sf mod} $h$; and $x \not\in L_1 \cup L_2$ iff 
$|(x)(R_1 \ {\sf or} \ R_2)| \equiv 0$ {\sf mod} $h$. Hence, 
$L_1, L_2 \in \oplus_{h,1,0} \! \bullet \! {\sf NP}$ \, implies \, 
$L_1 \cup L_2 \in \oplus_{h,\{1,2\},0} \! \bullet \! {\sf NP}$.

Replacing $L_1, L_2$ by their complements $\overline{L}_1, \overline{L}_2$, 
and using the fact that $\oplus_{h,1,0} \! \bullet \! {\sf NP}$ is closed 
under complement, we also have: 
 \  $L_1, L_2 \in \oplus_{h,1,0} \! \bullet \! {\sf NP}$
 \, implies \,  $\overline{L}_1, \overline{L}_2 \in$
$ \oplus_{h,1,0} \! \bullet \! {\sf NP}$, which by the above implies
$\overline{L_1 \cap L_2} \in \oplus_{h,\{1,2\},0} \! \bullet \! {\sf NP}$.
Since ${\sf co} \oplus_{h,\{1,2\},0} \! \bullet \! {\sf NP} = $
$\oplus_{h,0,\{1,2\}} \! \bullet \! {\sf NP}$, we have:
$L_1 \cap L_2 \in \oplus_{h,0,\{1,2\}} \! \bullet \! {\sf NP}$. 
By Lemma \ref{h_i_Add}, $\oplus_{h,0,\{1,2\}} \! \bullet \! {\sf NP} = $
$\oplus_{h,1,\{2,3\}} \! \bullet \! {\sf NP}$. Also, 
$\oplus_{h,1,\{2,3\}} \! \bullet \! {\sf NP} \subseteq $
$\oplus_{h,1} \! \bullet \! {\sf NP}$ ($ = \oplus_h \! \bullet \! {\sf NP}$). 
Hence, $L_1 \cap L_2 \in \oplus_h \! \bullet \! {\sf NP}$.

Since we saw in Lemma \ref{NPcoNPsubPlusNP} that {\sf NP} and {\sf coNP} 
are contained in $\oplus_{h,1,0} \! \bullet \! {\sf NP}$, it follows 
that ${\sf DP}$ is contained in $\oplus_h \! \bullet \! {\sf NP}$.   
  \ \ \ $\Box$

\medskip

We will see next that $\oplus_h \! \bullet \! {\sf NP}$ and 
$\oplus_{h,1,0} \! \bullet \! {\sf NP}$ have complete problems (with respect 
to polynomial-time many-to-one reduction). On the other hand {\sf BH} and 
{\sf PH} do not have complete problems -- unless the these hierarchies 
collapse. 
It is also known that a collapse of {\sf BH} implies a collapse of 
{\sf PH} (Kadin and Chang \cite{Kadin,ChangKadin}). 
This shows that we have:
 \ {\it Unless the polynomial hierarchy {\sf PH} collapses, 
$\oplus_h \! \bullet \! {\sf NP}$ and $\oplus_{h,1,0} \! \bullet \! {\sf NP}$ 
are both different from {\sf BH} and different from {\sf PH}. 
}

Recall the $\forall \exists$-quantified boolean formula problem, also called
$\forall \exists${\sc Sat}. The input for $\forall \exists${\sc Sat} is a
fully quantified boolean formula 
 \, $(\forall y_1, \ldots, y_n) \, (\exists x_1, \ldots, x_m)$
$B(x_1, \ldots, x_m, y_1, \ldots, y_n)$,  
and the question is whether this formula is true. It is well known that 
$\forall \exists${\sc Sat} is $\Pi_2^{\sf P}$-complete. In a similar way,
{\sc Sat} can be extended by any other quantifier sequence, which provides
complete problems for the classes $\Pi_{\ell}^{\sf P}$ and 
$\Sigma_{\ell}^{\sf P}$ of the polynomial hierarchy {\sf PH}. Another 
extension of {\sc Sat}, called $\#${\sc Sat} (``number sat''), is complete 
in the class $\# {\sf P}$ for parsimonious many-to-one polynomial-time 
reductions (Valiant \cite{Valiant1,Valiant2}; see also Chapter 8 of 
\cite{Papadim} for the definition of these reductions). 
The problem $\#${\sc Sat} is the function which maps any boolean formula 
$B(x_1, \ldots, x_m)$ to
 \ $|\{ (b_1, \ldots, b_m) \in \{0,1\}^m : B(b_1, \ldots, b_m) = 1\}|$  
 \ (i.e., the number of satisfying truth-value assignments, this number
being represented in binary).
This was generalized by \cite{DurandHerKol} to $\#\Pi_{\ell}${\sc Sat} and 
$\#\Sigma_{\ell}${\sc Sat} which are complete in 
$\# \bullet \Pi_{\ell}^{\sf P}$, respectively 
$\# \bullet \Sigma_{\ell}^{\sf P}$ (again for parsimonious many-to-one 
polynomial-time reductions). 

In the context of $\oplus_h \! \bullet \! {\sf NP}$ we introduce the following 
extension of {\sc Sat}, called $\oplus_h \exists${\sc Sat}. \\  
{\sf Input:} \ An existentially quantified boolean formula 
$(\exists x_1, \ldots, x_m ) B(x_1, \ldots, x_m, y_1, \ldots, y_n)$ with
free variables $y_1, \ldots, y_n$, 
where $B(x_1, \ldots, x_m, y_1, \ldots, y_n)$ is an ordinary boolean formula 
whose variables range over $\{0,1\}$.   \\  
{\sf Question ($\oplus_h \exists${\sc Sat}-problem):} 
 \ Does the following hold: 
 
\smallskip 

$|\{ (b_1, \ldots, b_n) \in \{0,1\}^n \ : \ $
  $(\exists x_1, \ldots, x_m ) \, B(x_1, \ldots, x_m, b_1, \ldots, b_n)\}|$
  $ \ \equiv \ \ 1$ \ {\sf mod} $h$ ?

\smallskip

\noindent In a similar way we define the problem
$\oplus_{h,1,0} \exists${\sc Sat}. \\    
{\sf Input:}
 \ $(\exists x_1, \ldots, x_m ) \, B(x_1, \ldots, x_m, y_1, \ldots, y_n)$, 
as in $\oplus_h \exists${\sc Sat}. \\   
{\sf Question ($\oplus_{h,1,0} \exists${\sc Sat}-problem):} 
 \ Does the following hold: 

\smallskip

$|\{ (b_1, \ldots, b_n) \in \{0,1\}^n \ : \ $ 
  $(\exists x_1, \ldots, x_m ) \, B(x_1, \ldots, x_m, b_1, \ldots, b_n)\}|$
   $ \ \equiv \ \ 1$ \ {\sf mod} $h$ ,  

\smallskip

{\bf and} 

\smallskip

$|\{ (b_1, \ldots, b_n) \in \{0,1\}^n \ : \ $ 
   {\sf not}$(\exists x_1, \ldots, $
   $x_m ) \, B(x_1, \ldots, x_m, b_1, \ldots, b_n)\}| \ \equiv \ \ 0$ 
   \ {\sf mod} $h$ {\bf ?}

\smallskip

\noindent The same parsimonious many-to-one polynomial-time reductions that
prove completeness of $\#${\sc Sat} and of $\#\Pi_{\ell}${\sc Sat}
yield the following: 

{\em The problem $\oplus_h \exists${\sc Sat} is 
 \, $\oplus_h \! \bullet \! {\sf NP}$-complete, and the problem 
$\oplus_{h,1,0} \exists${\sc Sat} is 
 \, $\oplus_{h,1,0} \! \bullet \! {\sf NP}$-complete. }

%\bigskip

%%%
\subsection{The image size problem}

In the following we will use combinational circuits (i.e., acyclic digital
circuits, made from {\sf and}, {\sf or}, {\sf not}, and {\sf fork} gates).
We will need to generalize these circuits to {\bf partial combinational
circuits}, simply by allowing a one-wire gate {\sf id}$_1$ which maps the
boolean value 1 to 1, and is undefined on 0. We will denote ``undefined'' 
by $\bot$.  When a gate has $\bot$ on one (or more) of its input wires, 
its output will be $\bot$. We also add the following rule about partial 
outputs of a circuit:

{\bf Partial outputs rule:}
 \ {\it If one or more output wires of a circuit receive the undefined 
value $\bot$ then the entire output of the circuit is viewed as undefined.}

In other words, any string in $\{0,1,\bot\}^*$ containing at least one 
$\bot$ is equivalent to $\bot$; so, up to this equivalence,
$\{0,1,\bot\}^*$  is $\{0,1\}^* \cup \{ \bot \}$.
The inputs of a partial 
combinational circuit $C$ are the elements of $\{0,1\}^m$ for some $m$ 
(depending on $C$). For $x \in \{0,1\}^m$ the output belongs to 
$\{0,1\}^n \cup \{ \bot \}$ for some $n$ (depending on $C$), and is 
denoted by $C(x)$.
We denote the domain of $C$ by ${\sf Dom}(C)$; it consists of the bitstrings 
in $\{0,1\}^m$ on which the output is defined. Hence, 
${\sf Dom}(C) = \{ x \in \{0,1\}^m : C(x) \neq \bot\}$. 
The set of all outputs of $C$ (not counting $\bot$) is called the image 
of $C$ and is denoted by ${\sf Im}(C)$; so, 
${\sf Im}(C) = \{ C(x) \in \{0,1\}^n : x \in \{0,1\}^m \}$.

To get closer to the $\equiv_{\cal D}$-decision problem of $M_{k,1}$ over
$\Gamma \cup \tau$, we introduce the following problems, called the 
{\it image size problem for partial combinational circuits} and the 
{\it modular image size problem for partial combinational circuits}. \\  
{\sf Input:} \ A partial combinational circuit $C$. \\     
{\sf Output} {\bf (image size problem):} \ The binary representation of 
the number $|{\sf Im}(C)|$ (i.e., the number of non-$\bot$ outputs; the 
outcome $\bot$, if it occurs, is not counted as an output). \\    
{\sf Question ({\sf mod} $h$ image size problem, for fixed $h \geq 2$):} 
 \ $|{\sf Im}(C)| \equiv 1$ {\sf mod} $h$ ? 

\smallskip

Finally, in relation to the $\equiv_{\cal D}$-decision problem we introduce
the {\bf modular image size problem of $M_{k,1}$} over 
$\Gamma \cup \tau$: \\   
{\sf Input:} \ $\varphi \in M_{k,1}$, given by a word over 
$\Gamma \cup \tau$. \\  
{\sf Question:}
 \ $|{\sf imC}(\varphi)| \equiv 1$ {\sf mod} $k-1$ ?
 
\smallskip

The number $|{\sf imC}(\varphi)|$ depends on the right ideal homomorphism 
that is chosen to represent $\varphi$; however, $|{\sf imC}(\varphi)|$ 
{\sf mod} $k-1$ does not depend the choice of representative 
(Prop.\ 2.4 in \cite{BiThomMon}); i.e., $|{\sf imC}(\varphi)|$ {\sf mod} 
$k-1$ is an invariant of $\varphi$ as an element of $M_{k,1}$. 
We only consider the modular image size problem of $M_{k,1}$ when 
$k \geq 3$. Recall that $M_{2,1}$ has only one non-zero $\cal D$-class.

\begin{thm} \label{complexity_imagesize} . \\   
{\bf (1)} The image size problem for partial combinational circuits is
$\# \bullet {\sf NP}$-complete.  \\  
{\bf (2)} The {\sf mod} $h$ image size problem for partial combinational 
circuits (for $h \geq 2$) is $\oplus_h \! \bullet \! {\sf NP}$-complete.  \\
{\bf (3)} For $k \geq 3$, the modular image size problem of $M_{k,1}$ over 
$\Gamma \cup \tau$ is $\oplus_{k-1} \! \bullet \! {\sf NP}$-complete. 
\end{thm}
{\bf Proof.} {\bf (1)} \ To prove that the image size problem is in
$\# \bullet {\sf NP}$ we consider the predicate $R$ defined by

\smallskip

 \ \ \ $(C,y) \in R$ \ \ iff 
 \ \ $(\exists x \in {\sf Dom}(C)) (C(x) = y)$, 

\smallskip

\noindent where $C$ ranges over all partial combinational circuits.
Clearly, the membership problem of $R$ is in {\sf NP}, and $R$ is 
polynomially balanced (in fact, $|y| \leq |C|$ since the output ports of 
$C$ are counted in the size of $C$). 
Then we have \  $\{ y : (C,y) \in R\} = {\sf Im}(C)$, hence the function 
\ $C \mapsto |{\sf Im}(C)|$ \ is in $\# \bullet {\sf NP}$. 

To prove $\# \bullet {\sf NP}$-hardness we will reduce 
$\# \exists${\sc Sat} to the image size problem. 
Let $B(x_1,x_2)$ be a boolean formula where $x_1$ is a sequence of $m$ 
boolean variables, and $x_2$ is a sequence of $n$ boolean variables.
We map $B$ to a partial combinational circuit $C_{B,m,n}$ with partial 
input-output function defined by

\smallskip

 \ \ \ $(x_1, x_2) \ \longmapsto \ C_{B,m,n}(x_1, x_2) \ = \ $
  $\left\{ \begin{array}{ll} 
           x_2 \  & \mbox{if \ $B(x_1, x_2) = 1$,} \\   
           \bot \ & \mbox{if \ $B(x_1, x_2) = 0$.}
           \end{array} \right. $

\smallskip

\noindent From the formula for $B(x_1, x_2)$ one can easily construct
a partial combinational circuit for $C_{B,m,n}$. Moreover, 
${\sf Im}(C_{B,m,n}) = \{x_2 : (\exists x_1) B(x_1,x_2) \}$, hence the 
reduction is a parsimonious reduction from the function 
 \ $\big[ \, B \longmapsto |\{ x_2 : (\exists x_1) B(x_1,x_2)\}| \, \big]$ 
 \ to the function 
 \ $\big[ \, C_{B,m,n} \longmapsto |{\sf Im}(C_{B,m,n})| \, \big]$.

\smallskip

\noindent {\bf (2)} \ Membership in $\oplus_h \! \bullet \! {\sf NP}$ is
proved as in {\bf (1)}. 
The reduction in {\bf (1)} also yields a parsimonious reduction of 
$\oplus_h${\sc Sat} to the {\sf mod} $h$ image size problem. This shows
$\oplus_h \! \bullet \! {\sf NP}$-hardness.

\smallskip

\noindent {\bf (3)} \ To prove that the modular image size problem of 
$M_{k,1}$ is in $\oplus_{k-1} \! \bullet \! {\sf NP}$ we consider the predicate 
$R$ defined by 

\smallskip

 \ \ \ $(\varphi,y) \in R$ \ \ iff \ \ $y \in {\sf imC}(\varphi)$. 

\smallskip

\noindent
Here, $\varphi$ is expressed by a word over $\Gamma \cup \tau$, where each
$\tau_{i-1,i} \in \tau$ has length $|\tau_{i-1,i}| = i$. 

The predicate $R$ is in {\sf NP}; see Prop.\ 4.9 about the image 
membership problem in \cite{BiRL}. 
The predicate $R$ is also polynomially balanced. In fact, for 
$y \in {\sf imC}(\varphi)$ we have by Theorem 4.5(2) in \cite{BiThomMon}:  
$|y| \leq c \cdot |\varphi|_{\Gamma \cup \tau}$ (for some constant $c>0$),
since we have $|\tau_{i-1,i}| = i$. 

Hardness follows from {\bf (2)}, since partial combinational circuits are
special cases of elements of $M_{k,1}$ expressed over $\Gamma \cup \tau$. 
 \ \ \ $\Box$

\medskip

\noindent {\bf Remark.} The proof of {\bf (1)} above shows why {\it partial}
circuits were introduced: For an ordinary (total) circuit $C$ the image
size is never 0, whereas the set $\{x_2 : (\exists x_1) B(x_1,x_2)\}$ 
can be empty. So there is no parsimonious reduction from 
$\# \exists${\sc Sat} to the image-size problem of total circuits.    

\bigskip

For comparison, the {\it domain size problem} for partial combinational 
circuits (i.e., the function $C \longmapsto |\{ x : x \in {\sf Dom}(C)\}|$)
is \#{\sf P}-complete.
Indeed, we can map any boolean formula $B$ to a partial combinational 
circuit $C_B$ which (on input $x$) outputs $\bot$ when $B(x) = 0$,
and outputs 1 when $B(x) = 1$. Then the domain of the partial circuit $C_B$
consists of the satisfying truth values of $B$, so this is a parsimonious 
reduction of \#{\sc Sat} to the domain size problem. 
Moreover, the domain size problem is in \#{\sf P}. Indeed, the predicate
$\{(x,C) : x \in {\sf Dom}(C)\}$ (where $x \in \{0,1\}^*$ and $C$ ranges 
over partial combinational circuits) is in {\sf P} since a circuit can be
evaluated quickly on a given input.
 
For a fixed $h \geq 2$ and $0 \leq i \leq h-1$
we can also consider the {\it modular domain size problem} for partial 
combinational circuits; for a circuit $C$, the question is whether 
 \ $|{\sf Dom}(C)| \equiv 1$ {\sf mod} $h$.
As above one proves that this problem is $\oplus_h {\sf P}$-complete. 

Similarly, for $1 \leq i \leq k-1$ we have the {\it modular domain code 
size problem} in $M_{k,1}$; for $\varphi \in M_{k,1}$, given by a word 
over $\Gamma \cup \tau$, the question is whether
 \ $|{\sf domC}(\varphi)| \equiv i$ {\sf mod} $k-1$.
This problem is $\oplus_{k-1} {\sf P}$-complete.

%%%
\subsection{The complexity of $\equiv_{\cal D}$ over $\Gamma \cup \tau$}

Recall that $M_{k,1}$ has $k-1$ non-zero $\cal D$-classes 
$D_i = \{ \varphi : |{\sf imC}(\varphi)| \equiv i$ {\sf mod} $k-1 \}$, for 
$1 \leq i \leq k-1$.

\begin{thm} \label{complexity_D_relation} \  
Let $k \geq 3$ and $1 \leq i \leq k-1$. The membership problem of the 
$\cal D$-class $D_i$ of $M_{k,1}$ over $\Gamma \cup \tau$ is \, 
$\oplus_{k-1} \! \bullet \! {\sf NP}$-complete.
\end{thm}
{\bf Proof.} Checking whether an element is not $\equiv_{\cal D} {\bf 0}$ 
is in {\sf NP} (by Prop.\  \ref{0_word_probl_infgen}), and  {\sf NP} is 
contained in $\oplus_{k-1} \! \bullet \! {\sf NP}$.
Checking whether a non-zero element is in $D_i$ is
$\oplus_{k-1} \! \bullet \! {\sf NP}$-complete by Theorem 
\ref{complexity_imagesize}(3), and by the fact that for a non-zero element
$\varphi \in M_{k,1}$ we have \ $\varphi \in D_i$
 \ iff \ $|{\sf imC}(\varphi)| \equiv i$ {\sf mod} $k-1$ 
 \ (Theorem 2.5 in \cite{BiThomMon}). 
 \ \ \ $\Box$

\medskip

\noindent {\bf Remark.} {\it The $\equiv_{\cal D}$-decision problem of
$M_{2,1}$ over $\Gamma \cup \tau$ is {\sf coDP}-complete, with respect to 
polynomial-time disjunctive reduction. } \ Indeed, in $M_{2,1}$, 
$\equiv_{\cal D}$ and $\equiv_{\cal J}$ are the same (Theorem 2.5 in 
\cite{BiThomMon}), and we saw in Prop.\ \ref{J_order_infgen} that the 
$\equiv_{\cal J}$-decision problem is {\sf coDP}-complete.

\bigskip

Earlier we considered the $\cal D$-{\it pivot search problem} of $M_{k,1}$, 
and we proved that it is in {\sf P} when inputs are expressed over a finite 
generating set $\Gamma$ of $M_{k,1}$.
Over circuit-like generating sets $\Gamma \cup \tau$ we have the following.

\begin{thm} \label{pivotLength_GammaTau} \   
The $\cal D$-pivots of $M_{3,1}$ do not have polynomially bounded 
word-length over $\Gamma \cup \tau$, unless the polynomial hierarchy 
{\sf PH} collapses.
More precisely, suppose there is a polynomial $p(.)$ such that for all 
$\psi, \varphi \in M_{3,1}$ we have: $\psi \equiv_{\cal D} \varphi$ implies
that there is a $\cal D$-pivot $\chi$ with 
 \ $|\chi|_{\Gamma \cup \tau} \ \leq \ $
$p(|\psi|_{\Gamma \cup \tau} + |\varphi|_{\Gamma \cup \tau})$; 
then {\sf PH} collapses to $\Pi_4^{\sf P} \cap \Sigma_4^{\sf P}$.
\end{thm}
{\bf Proof.} We proved in \cite{BiRL} that the $\equiv_{\cal R}$- and 
$\equiv_{\cal L}$-decision problems are in $\Pi_2^{\sf P}$. 
If $\cal D$-pivots had polynomially bounded word-length over 
$\Gamma \cup \tau$ then the $\equiv_{\cal D}$-decision problem would 
be in $\Sigma_3^{\sf P}$, by just guessing a pivot $\chi$ in 
nondeterministic polynomial time, and checking whether 
$\psi \equiv_{\cal L} \chi \equiv_{\cal R} \varphi$ (which is in 
$\Pi_2^{\sf P}$).
However, the $\equiv_{\cal D}$-decision problem is complete in 
$\oplus_{2} \! \bullet \! {\sf NP}$, hence 
$\oplus_{2} \! \bullet \! {\sf NP}$ would be in $\Sigma_3^{\sf P}$; 
hence  $\oplus {\sf P}$ would be contained in $\Sigma_3^{\sf P}$. 
By \cite{Toda91}, $\oplus {\sf P} \subseteq \Sigma_3^{\sf P}$ implies that 
{\sf PH} collapses to $\Pi_4^{\sf P} \cap \Sigma_4^{\sf P}$.
 \ \ \ $\Box$

\bigskip

For $M_{k,1}$ with $k > 3$, Theorem \ref{pivotLength_GammaTau} probably 
also holds, but the {\sf mod} $h$ version of Toda's theorem (for $h>2$)
has not been checked. 
% \ {\sl [Redo the proof in pp.\ 334-340 of \cite{DuKo} for $\oplus_h {\sf P}$;
%check if it works. An alphabet of size $h$ is now used (instead of $\{0,1\}$),
%in order to define the operators $\oplus_h$, {\it BP}, and the classes
%$\oplus_h {\sf P}$ and {\it BPP}.]   }

Case $k = 2:$ We leave it as an open question whether Theorem 
\ref{pivotLength_GammaTau} holds for $M_{2,1}$.
All non-zero elements of $M_{2,1}$ are $\cal D$-equivalent, so here pivots
always exist.

%%%%%%%%%%%%%%%%%%%%%%%%%%%%%%%%%%%%%%%%%%%%%%%%%%%
%% Section
%%%%%%%%%%%%%%%%%%%%%%%%%%%%%%%%%%%%%%%%%%%%%%%%%%%

\section{The Green relations of ${\it Inv}_{k,1}$ and their complexity }

We saw in \cite{BiThomMon} that ${\it Inv}_{k,1}$ is an {\it inverse} 
monoid, i.e., for every $\alpha \in M$ there exists one and only one 
$\alpha' \in M$ such that $\alpha \alpha' \alpha = \alpha$ and 
$\alpha' \alpha \alpha' = \alpha'$; the element $\alpha'$ is called the
inverse of $\alpha$. Some elementary facts about inverse monoids:  
 \ For all $\alpha, \beta  \in M:$ 
$(\alpha \cdot \beta)' = \beta' \cdot \alpha'$. 
For all $\alpha, \beta  \in M:$ $\beta \leq_{\cal R} \alpha$ iff
$\alpha' \leq_{\cal L} \beta'$; similarly, $\beta \leq_{\cal L} \alpha$ iff
$\alpha' \leq_{\cal R} \beta'$.

\begin{pro} \label{GreenRelInvk1} {\bf (The Green relations of 
  ${\it Inv}_{k,1}$).}  \\  
{\bf (1)} \ The Green relations $\leq_{\cal J}, \equiv_{\cal D}$,
$\leq_{\cal R}, \leq_{\cal L}$ of ${\it Inv}_{k,1}$ are the restrictions
of the corresponding Green relations of $M_{k,1}$.   \\   
{\bf (2)} \ ${\it Inv}_{k,1}$ is a union of $\equiv_{\cal L}$-classes 
of $M_{k,1}$.  In other words, if an $\cal L$-class of $M_{k,1}$ intersects
${\it Inv}_{k,1}$ then this entire $\cal L$-class is contained in 
${\it Inv}_{k,1}$.    \\   
{\bf (3)} \ Every $\equiv_{\cal R}$-class of $M_{k,1}$ intersects 
${\it Inv}_{k,1}$. \\     
{\bf (4)} \ Let $\Gamma_I$ be a finite generating set of ${\it Inv}_{k,1}$,
and let us assume that $\Gamma_I$ is closed under inverse. 
Then for every $\varphi \in {\it Inv}_{k,1}$ we have 
 \ $|\varphi|_{\Gamma_I} = |\varphi^{-1}|_{\Gamma_I}$ \ and 
 \ $|\varphi|_{\Gamma_I \cup \tau} = |\varphi^{-1}|_{\Gamma_I \cup \tau}$.  
\end{pro}
{\bf Proof.} Let us use ${\cal L}(M_{k,1})$ to indicate the 
$\cal L$ relations of $M_{k,1}$, and similarly for ${\cal R}(M_{k,1})$. 

\smallskip

\noindent {\bf (1)} \ This is Lemma 2.9 in \cite{BiThomMon}. 

\smallskip

\noindent {\bf (2)} \ For $\varphi \in M_{k,1}$ we have: \    
$\varphi \in {\it Inv}_{k,1}$ iff ${\sf part}(\varphi)$ is the identity
congruence on ${\sf Dom}(\varphi)$. By the characterization of 
$\leq_{{\cal L}(M_{k,1})}$ in $M_{k,1}$ (Theorem 3.32 in the arXiv version 
of \cite{BiRL}), this implies that every element in the 
${\cal L}(M_{k,1})$-class of $\varphi$ has the identity congruence as its 
partition.  Hence the whole ${\cal L}(M_{k,1})$-class is contained in 
${\it Inv}_{k,1}$.

\smallskip

\noindent {\bf (3)} \ Let $\varphi: P \to Q$ be a table for an element of
$M_{k,1}$, where $P$ and $Q$ are finite prefix codes. Then ${\sf id}_Q$ 
belongs to the ${\cal R}(M_{k,1})$-class of $\varphi$ (by Theorem 2.1 in 
\cite{BiRL}), and ${\sf id}_Q$ obviously belongs to ${\it Inv}_{k,1}$. 

\smallskip

\noindent {\bf (4)} \ This is straightforward. 
 \ \ \ $\Box$

\begin{pro} \label{complexityGrRelInvFG} \ The decision problems for the 
Green relations $\leq_{\cal J}, \equiv_{\cal D}$,
$\leq_{\cal R}, \leq_{\cal L}$ of ${\it Inv}_{k,1}$ are in {\sf P} 
when the inputs are given by words over a finite generating set of 
${\it Inv}_{k,1}$.
\end{pro}
{\bf Proof.} Since ${\it Inv}_{k,1}$ is a finitely generated submonoid of 
$M_{k,1}$, this is a consequence of the corresponding result for $M_{k,1}$
(Proposition \ref{J_D_order_probl_fg} above and Theorems 5.1 and 6.1 in
\cite{BiRL}). 
 \ \ \ $\Box$

\medskip

As a consequence of Prop.\ \ref{GreenRelInvk1}(4) and the elementary facts
about inverse monoids mentioned before Prop.\ \ref{GreenRelInvk1}, the 
$\leq_{\cal R}$ decision problem and the $\leq_{\cal L}$ decison problem of 
${\it Inv}_{k,1}$ (over a circuit-like alphabet $\Gamma_I \cup \tau$) can 
be reduced to each other and have the same computational complexity. 

\medskip

Let $\Gamma_I$ be a finite generating set of ${\it Inv}_{k,1}$; we can 
assume that $\Gamma_I$ is closed under inverse (since this is only a 
finite change in the generating set).
The {\it {\bf 0} word problem} of ${\it Inv}_{k,1}$ over 
$\Gamma_I \cup \tau$ is specified as follows.

\smallskip

\noindent {\sf Input:} \ $\varphi \in {\it Inv}_{k,1}$, given by a word 
over the generating set $\Gamma_I \cup \tau$,

\smallskip

\noindent {\sf Question}:
 \ Is \ $\varphi = {\bf 0}$ \ as an element of ${\it Inv}_{k,1}$ ?

\begin{thm} \label{zeroWPforInv} \  
The {\bf 0} word problem of ${\it Inv}_{k,1}$ over $\Gamma_I \cup \tau$ 
is {\sf coNP}-complete.
\end{thm}
{\bf Proof.} \ The problem is in {\sf coNP}, for the same reason as the
{\bf 0} word problem of $M_{k,1}$ over $\Gamma \cup \tau$ is in {\sf coNP} 
(Prop.\ 6.2 in \cite{BiRL}). 

To show {\sf coNP}-hardness we reduce the tautology problem for boolean
formulas to the {\bf 0} word problem of ${\it Inv}_{k,1}$. This is done in 
two steps; in the first step we work over the alphabet $\{0,1\}$, rather 
than over $A = \{a_1, \ldots, a_k\}$, and (if $k>2$) in the second step we 
use $A$. 

Let $\Gamma_{I,k}$ and $\Gamma_{G,k}$ be finite generating sets for, 
respectively, ${\it Inv}_{k,1}$ and $G_{k,1}$. We can assume that 
$\Gamma_{G,k} \subset \Gamma_{I,k}$ (since only finite changes are needed 
to achieve this).

\smallskip

\noindent {\sf Step 1.} \ We will reduce the tautology problem for boolean
formulas to the {\bf 0} word problem of ${\it Inv}_{2,1}$ (over 
$\Gamma_{I,2} \cup \tau$). Let $B(x_1, \ldots, x_m)$ be any boolean 
formula; it defines a map $B: \{0,1\}^m \to \{0,1\}$.
By Theorem 4.1 in \cite{BiDistor}, we map $B$ (given by a boolean 
formula or a circuit) to an element $\Phi_B \in G_{2,1}$ (given by a word 
over $\Gamma_{G,2} \cup \tau$), such that for all $x \in \{0,1\}^m$:

\smallskip

 \ \ \ $\Phi_B(0 x) \ = \ 0 \ B(x) \ x$.

\smallskip

\noindent By Theorem 4.1 in \cite{BiDistor}, this mapping from a formula
$B$ to word for $\Phi_B$ can be computed in deterministic polynomial time.  
Now we have: 

\smallskip

 \ \ \ $B$ is a tautology \ \ iff 
 \ \ ${\sf id}_{0 0 \{0,1\}^*} \circ \Phi_B \circ {\sf id}_{0\{0,1\}^*} $
     $ \ = \ {\bf 0}$. 

\smallskip

\noindent Since $\Gamma_{G,2} \subset \Gamma_{I,2}$, $\Phi_B$ is 
automatically over $\Gamma_{I,2} \cup \tau$.
Also, the partial identities ${\sf id}_{0\{0,1\}^*}$ and
${\sf id}_{0 0 \{0,1\}^*}$ are fixed elements of ${\it Inv}_{2,1}$ and 
they can be represented by fixed words over $\Gamma_{I,2}$. So the map 
 \ $B \ \longmapsto \ $
   ${\sf id}_{0 0 \{0,1\}^*} \circ \Phi_B \circ {\sf id}_{0\{0,1\}^*}$  
 \ reduces the tautology problem for boolean formulas to the {\bf 0} word 
problem of ${\it Inv}_{2,1}$ (over $\Gamma_{I,2} \cup \tau$).

\smallskip

\noindent {\sf Step 2.} \   
We reduce the {\bf 0} word problem of ${\it Inv}_{2,1}$ (over
$\Gamma_{I,2} \cup \tau$) to the {\bf 0} word problem of ${\it Inv}_{k,1}$ 
(over $\Gamma_{I,k} \cup \tau$), for any $k \geq 2$.
Let $\psi \in {\it Inv}_{2,1}$, and let \ $\ell(\psi) = $
${\sf max}\{ |z| : z \in {\sf domC}(\psi) \cup {\sf imC}(\psi)\}$.
Suppose $\psi$ is given by a word $w$ over $\Gamma_{I,2} \cup \tau$.

Let $\gamma \in \Gamma_{I,2}$ be any generator, and let us take a table
$P \to Q$ be for $\gamma$, where $P, Q \subset \{0,1\}^*$ are finite 
prefix codes. 
We view $\gamma$ as an element $\gamma^A$ of ${\it Inv}_{k,1}$ by taking the 
table $P \to Q$ as a table over $A = \{a_1, \ldots, a_k\}$, by 
identifying $\{0,1\}$ with $\{a_1, a_2\} \subseteq A$. Let 
$\Gamma^A_{I,2} = \{\gamma^A : \gamma \in \Gamma_{I,2}\}$. 
Since $\Gamma^A_{I,2}$ is finite we can assume that 
$\Gamma^A_{I,2} \subset \Gamma_{I,k}$ (since only finite changes are needed
to achieve this).

Let $W$ be the word over $\Gamma_{I,2} \cup \tau$ obtained by replacing 
each generator $\gamma \in \Gamma_{I,2}$ by the corresponding $\gamma^A$; 
elements of $\tau$ are not changed (except that they now act on $A^*$, 
rather than just $\{0,1\}^*$).
Let $\Psi \in {\it Inv}_{k,1}$ be the element of ${\it Inv}_{k,1}$ 
represented by $W$. For $\psi \in {\it Inv}_{2,1}$, \, $\psi(z)$ is 
undefined when $z \not\in \{a_1, a_2\}^*$.  
For a prefix code $P \subset A^*$ we abbreviate ${\sf id}_{P A^*}$ to
${\sf id}_P$.    Then we have:

\smallskip

\noindent {\sf Claim.} \ \ For all $z \in A^{\ell(\psi)} :$ 
 \ $\psi(z) = \Psi \circ {\sf id}_{\{a_1,a_2\}^{\ell(\psi)}}(z)$.

 \ \ \ \ \ \ Moreover, $\psi = {\bf 0}$ as an element of ${\it Inv}_{2,1}$ 
 \ iff \ $\Psi \circ {\sf id}_{\{a_1,a_2\}^{\ell(\psi)}} = {\bf 0}$ 
 \ as an element of ${\it Inv}_{k,1}$.

\smallskip

\noindent Indeed, both sides of the equality are undefined on 
$A^{\ell(\psi)} - \{a_1,a_2\}^{\ell(\psi)}$. 
For $z \in \{a_1,a_2\}^{\ell(\psi)}$ $( = \{0,1\}^{\ell(\psi)} )$ we have:
 \ ${\sf id}_{\{a_1,a_2\}^{\ell(\psi)}}(z) = z$ \ and 
 \ $\Psi(z) = \psi(z)$.  Moreover, both ${\sf domC}(\psi)$ and 
${\sf domC}(\Psi \circ {\sf id}_{ \{a_1,a_2\}^{\ell(\psi)} })$ are subsets 
of $A^{\leq \ell(\psi)}$. 
It follows that $\psi = {\bf 0}$ in ${\it Inv}_{2,1}$ \   iff 
 \ $\Psi \circ {\sf id}_{\{a_1,a_2\}^{\ell(\psi)}} = {\bf 0}$ 
 \ in ${\it Inv}_{2,1}$.
[This proves the Claim.]

\smallskip

\noindent One easily verifies that  

\smallskip

${\sf id}_{\{a_1,a_2\}^{\ell(\psi)}} \ = \ $
$\tau_{\ell(\psi),1} \circ {\sf id}_{\{a_1,a_2\}} \circ $
$\tau_{\ell(\psi),1} \circ \ \ldots \ \circ $
$\tau_{j,1} \circ {\sf id}_{\{a_1,a_2\}} \circ \tau_{j,1} \circ $ 
$ \ \ldots \ \circ $
$\tau_{2,1} \circ {\sf id}_{\{a_1,a_2\}} \circ \tau_{2,1}(.)$.

\smallskip

\noindent 
Hence, the word-length of ${\sf id}_{\{a_1,a_2\}^{\ell(\psi)}}$ over
$\Gamma_{I,k} \cup \tau$ is polynomially bounded (in terms of
$|\psi|_{\Gamma_{I,2} \cup \tau}$).  

Thus the map from $\psi \in {\it Inv}_{2,1}$ (given by a word over
$\Gamma_{I,2} \cup \tau$) to $\Psi \circ {\sf id}_{\{a_1,a_2\}^{\ell(\psi)}}$ 
$\in {\it Inv}_{k,1}$ (given by a word over $\Gamma_{I,k} \cup \tau$) is 
polynomial-time computable.    Hence this map is a polynomial-time 
reduction from the {\bf 0} word problem of 
${\it Inv}_{2,1}$ to the {\bf 0} word problem of ${\it Inv}_{k,1}$.
 \ \ \ $\Box$

\begin{thm} \label{complexityLRinInv} {\bf (The $\cal R$ and $\cal L$
decision problems).}  \   
The $\leq_{\cal L}$ and $\leq_{\cal R}$ decision problems of 
${\it Inv}_{k,1}$ over $\Gamma_I \cup \tau$ are each {\sf coNP}-complete.
\end{thm}
{\bf Proof.} \ The $\leq_{\cal L}$ decision problem is in {\sf coNP} 
for $M_{k,1}$ (by Theorem 6.7 in \cite{BiRL}), hence (by Prop.\
\ref{GreenRelInvk1}(1)) it is in {\sf coNP} for ${\it Inv}_{k,1}$ too. 
In ${\it Inv}_{k,1}$, the $\leq_{\cal R}$ decision problem reduces to the 
$\leq_{\cal L}$ decision problem by Prop.\ \ref{GreenRelInvk1}(4), so the
$\leq_{\cal R}$ decision problem of ${\it Inv}_{k,1}$ is in {\sf coNP}.

The $\leq_{\cal L}$ and the $\leq_{\cal R}$ decision problems are 
{\sf coNP}-hard by Theorem \ref{zeroWPforInv}, since 
$\varphi \leq_{\cal L} {\bf 0}$ iff $\varphi = {\bf 0}$ (and similarly for
$\leq_{\cal R}$). 
 \ \ \ $\Box$

\begin{thm} \label{complexityJinInv} {\bf (The $\cal J$ decision problem).}  
 \ The $\leq_{\cal J}$ and the $\equiv_{\cal J}$ decision problems of 
${\it Inv}_{k,1}$ over $\Gamma_I \cup \tau$ are {\sf coDP}-complete (for the
$\equiv_{\cal J}$ decision problems the {\sf coDP}-completeness is with
respect to polynomial-time disjunctive reduction).
\end{thm}
{\bf Proof.} \ The proof is the same as for $M_{k,1}$ (Theorem 
\ref{J_order_infgen}).  For $\varphi, \psi \in {\it Inv}_{k,1}$ we have
$\psi \leq_{\cal J} \varphi$ iff $\psi = {\bf 0}$ or $\varphi \neq {\bf 0}$.
The result follows since the {\bf 0} word problem is {\sf coNP}-complete.
Moreover, the $\leq_{\cal J}$ decision problem reduces to the $\equiv_{\cal J}$
decision problem by a polynomial-time disjunctive reduction, since 
$\psi \leq_{\cal J} \varphi$ iff $\psi \equiv_{\cal J} {\bf 0}$ or 
$\psi \equiv_{\cal J} \varphi$.
 \ \ \ $\Box$

\medskip

We will prove next that the membership problem of a non-zero $\cal D$-class 
is easier for ${\it Inv}_{k,1}$ than for $M_{k,1}$ (seen in Theorem 
\ref{complexity_D_relation}), if $\oplus_{k-1}${\sf P} is different
from $\oplus_{k-1} \! \bullet \! {\sf NP}$. The class $\oplus_{h,i} {\sf P}$
(for integers $h \geq 2$ and $0 \leq i \leq h-1$) was defined at the 
beginning of Section 6.1. Just as for $\oplus_{h,i} \! \bullet \! {\sf NP}$,
we can prove that $\oplus_{h,i} {\sf P} = \oplus_{h,j} {\sf P}$ for all
$i,j$; therefore we denote $\oplus_{h,i} {\sf P}$ by $\oplus_h {\sf P}$ 
for every $i$.

\begin{thm} \label{complexityDinInv} {\bf (Complexity of $\cal D$).}
 \ Let $k \geq 3$ and $1 \leq i \leq k-1$. The membership problem of the
$\cal D$-class $D_i$ of ${\it Inv}_{k,1}$ over $\Gamma_I \cup \tau$ is \,
$\oplus_{k-1} {\sf P}$-complete.
\end{thm}
{\bf Proof.} \ Let $\varphi \in {\it Inv}_{k,1}$ be given by a word over
$\Gamma_I \cup \tau$. By injectiveness,
 \ $|{\sf imC}(\varphi)| = |{\sf domC}(\varphi)|$.   
Hence, by the characterization of the $\cal D$ relation (Theorem 2.5 in
\cite{BiThomMon}) and by the fact that this characterization applies to 
${\it Inv}_{k,1}$ as well (Prop.\ \ref{GreenRelInvk1}(1)), the 
$\cal D$-class $D_i$ of ${\it Inv}_{k,1}$ satisfies

\smallskip

$D_i \ = \ \{ \varphi \in {\it Inv}_{k,1} : $
   $|{\sf domC}(\varphi)| \equiv i$ {\sf mod} $k-1 \}$.

\smallskip

\noindent Since the predicate 
 \ $R = \{(x, \varphi) \in A^* \times (\Gamma_I \cup \tau)^* : $
      $ x \in {\sf domC}(\varphi) \}$ \ is in {\sf P} (by Prop.\ 5.6(1)
in \cite{BiRL}), we conclude that the membership problem of $D_i$ is
in $\oplus_{k-1} {\sf P}$.

To show that the membership problem of $D_i$ is $\oplus_{k-1} {\sf P}$-hard
we will reduce $\oplus_{k-1}${\sc Sat} to it by a polynomial-time 
parsimonious reduction. The input to $\oplus_{k-1}${\sc Sat} is any 
boolean formula $B(x_1, \ldots, x_m)$; this formula defines a boolean 
function $B: \{0,1\}^m \to \{0,1\}$. 

Let $\Gamma_{I,k}$ and $\Gamma_{G,k}$ be finite generating sets for, 
respectively, ${\it Inv}_{k,1}$ and $G_{k,1}$. We can
assume that $\Gamma_{G,k} \subset \Gamma_{I,k}$
(since this can be achieved by finite changes).
We build the reduction in two steps, the first for $k=2$, the second for
$k > 2$.

\smallskip

\noindent {\sf Step 1.} \  
As in the proof of Theorem \ref{zeroWPforInv}, we map $B$ to the element 
$\Phi_B \in G_{2,1}$ (given by a word
over $\Gamma_{G,2} \cup \tau$) such that for all $x \in \{0,1\}^m$:

\smallskip

 \ \ \ $\Phi_B(0 x) \ = \ 0 \ B(x) \ x$.

\smallskip

\noindent By Theorem 4.1 in \cite{BiDistor}, this mapping from a formula
$B$ to word for $\Phi_B$ can be computed in deterministic polynomial time.
Now we consider the element $\varphi_B \in {\it Inv}_{2,1}$ (given by a word
over $\Gamma_{I,2} \cup \tau$), defined by

\smallskip

 \ \ \ $\varphi_B(.) \ = \ $
${\sf id}_{0 1 \{0,1\}^*} \circ \Phi_B \circ {\sf id}_{0 \{0,1\}^*}(.)$

\smallskip

\noindent Then, 

\smallskip

 \ \ \ ${\sf domC}(\varphi_B) \ = \ 0 \, \{ x \in \{0,1\}^m : B(x)=1\}$ 
$ \ \subseteq \ \{0,1\}^{m+1}$ \ \ and 

\smallskip

 \ \ \ ${\sf imC}(\varphi_B) \ = \ 0 \, 1 \, \{x \in \{0,1\}^m : B(x)=1\}$
$ \ \subseteq \ \{0,1\}^{m+2}$.  

\smallskip

\noindent Hence, 

\smallskip

 \ \ \ $|{\sf imC}(\varphi_B)| \ = \ |\{x \in \{0,1\}^m : B(x) = 1\}|$.

\smallskip

\noindent Moreover, the partial identities ${\sf id}_{0 1 \{0,1\}^*}$ and 
${\sf id}_{0 \{0,1\}^*}$ can be given by fixed words over $\Gamma_{I,2}$.
So the map which sends a formula for $B$ to a word that represents 
$\varphi_B$ (over $\Gamma_{I,2} \cup \tau$) is polynomial-time computable, 
and it is parsimonious (in the sense that the image code size of $\varphi_B$ 
is the number of satisfying truth-value assignments of $B$).

\smallskip

\noindent {\sf Step 2.} \   
We identify $\{0,1\}$ with $\{a_1,a_2\} \subset A$. We will map a word
representing $\varphi_B$ (over $\Gamma_{I,2} \cup \tau$) to a word over
$\Gamma_{I,k} \cup \tau$, representing an element 
$\psi_B \in {\it Inv}_{k,1}$; this map should be polynomial-time computable,
and it should be parsimonious in the sense that 
 \ $|{\sf imC}(\varphi_B)| \ = \ |{\sf imC}(\psi_B)|$.

Let $w \in (\Gamma_{I,2} \cup \tau)^*$ be a word that represents 
$\varphi_B$. Let $\gamma \in \Gamma_{I,2}$ be any generator, and let 
$P \to Q$ be a table for $\gamma$, where $P, Q \subset \{0,1\}^*$ are finite 
prefix codes.
We view $\gamma$ as an element $\gamma^A$ of ${\it Inv}_{k,1}$ by taking the
table $P \to Q$ as a table over $A = \{a_1, \ldots, a_k\}$, by identifying 
$\{0,1\}$ with $\{a_1, a_2\} \subseteq A$. Let
$\Gamma^A_{I,2} = \{\gamma^A : \gamma \in \Gamma_{I,2}\}$.
Since $\Gamma^A_{I,2}$ is finite we can assume that
$\Gamma^A_{I,2} \subset \Gamma_{I,k}$.

Let $W$ be the word over $\Gamma_{I,2} \cup \tau$ obtained from $w$ by
replacing each generator $\gamma \in \Gamma_{I,2}$ by the corresponding
$\gamma^A$; elements of $\tau$ are not changed (except that they now act on 
$A^*$).  
Let $\Phi_B$ be the element of ${\it Inv}_{k,1}$ represented by $W$. 
For $\varphi_B \in {\it Inv}_{k,1}$, \, $\varphi_B(z)$ is undefined when
$z \not\in \{a_1, a_2\}^*$.    We have:

\smallskip

\noindent {\sf Claim.} \ For all $z \in A^{m+1} :$
 \ $\varphi_B(z) = \Phi_B \circ {\sf id}_{\{a_1,a_2\}^{m+1} }(z)$.
 \ Moreover, after a restriction,   

\smallskip

 \ \ \ \ \  ${\sf imC}(\varphi_B) \ = \ $ 
${\sf imC}(\Phi_B \circ {\sf id}_{\{a_1,a_2\}^{m+1}}) \ $
$\subseteq \ a_1 a_2 \{a_1,a_2\}^m \ \subseteq \ \{a_1,a_2\}^{m+2}$. 

\smallskip

\noindent Proof of the Claim: \ By the definition of $\varphi_B$ we have
 \ ${\sf domC}(\varphi_B) \subseteq a_1 \{a_1,a_2\}^m$. 
Both sides of the equality are undefined on $A^{m+1} - \{a_1,a_2\}^{m+1}$.
For $z \in \{a_1,a_2\}^{m+1}$ we have:
 \ ${\sf id}_{\{a_1,a_2\}^{m+1}}(z) = z$ \ and
 \ $\Phi_B(z) = \varphi_B(z)$ \ (since $\Phi_B$ and $\varphi_B$ agree on
$\{a_1,a_2\}^*$). 
 \ \ \    [This proves the Claim.]

\smallskip

\noindent One easily verifies that

\smallskip

${\sf id}_{\{a_1,a_2\}^{m+1}} \ = \ $
$\tau_{m+1,1} \circ {\sf id}_{\{a_1,a_2\}} \circ \tau_{m+1,1} \circ $
$ \ \ldots \ \circ $
$\tau_{j,1} \circ {\sf id}_{\{a_1,a_2\}} \circ \tau_{j,1} \circ $
$ \ \ldots \ \circ $
$\tau_{2,1} \circ {\sf id}_{\{a_1,a_2\}} \circ \tau_{2,1}(.)$.

\smallskip

\noindent
Hence, the word-length of ${\sf id}_{\{a_1,a_2\}^{m+1}}$ over
$\Gamma_{I,k} \cup \tau$ is polynomially bounded.
Thus the map from $\varphi_B \in {\it Inv}_{2,1}$ (given by a word over
$\Gamma_{I,2} \cup \tau$) to \ $\Phi_B \circ {\sf id}_{\{a_1,a_2\}^{m+1}}$
$\in {\it Inv}_{k,1}$ (given by a word over $\Gamma_{I,k} \cup \tau$) is
polynomial-time computable.

Since it follows from the Claim that \ $|{\sf imC}(\varphi_B)| \ = \ $
$|{\sf imC}(\Phi_B \circ {\sf id}_{\{a_1,a_2\}^{m+1}})|$, the map from
$\varphi_B$ to \ $\Phi_B \circ {\sf id}_{\{a_1,a_2\}^{m+1}}$ \ is 
parsimonious.

Combining Step 1 and Step 2, we obtain a polynomial-time reduction from
the a boolean formula $B$ to an element 
$\Phi_B \circ {\sf id}_{\{a_1,a_2\}^{m+1}} \in {\it Inv}_{k,1}$ (given by a 
word over $\Gamma_{I,k} \cup \tau$). The reduction is parsimonious since
 \ $|\{x \in \{0,1\}^m : B(x)=1\}| \ = \ $
$|{\sf imC}(\Phi_B \circ {\sf id}_{\{a_1,a_2\}^{m+1}})|$.
Obviously, the latter equality holds modulo $k-1$ too.  
 \ \ \ $\Box$

\begin{thm} \label{pivotsizeInv} \  
The $\cal D$-pivots of ${\it Inv}_{3,1}$ do not have polynomially bounded
word-length over $\Gamma_I \cup \tau$, unless the polynomial hierarchy
{\sf PH} collapses to $\Sigma_3^{\sf P} \cap \Pi_3^{\sf P}$.
\end{thm}
{\bf Proof.} \ The proof is similar to the proof of Theorem 
\ref{pivotLength_GammaTau}. We saw that the $\equiv_{\cal L}$ and 
$\equiv_{\cal R}$ decision problems are in {\sf coNP}. If pivots always had
polynomially bounded lengths, the $\equiv_{\cal D}$ decision problem would
be in $\Sigma_2^{\sf P}$, by just guessing a pivot $\chi$ in
nondeterministic polynomial time, and checking whether
$\psi \equiv_{\cal L} \chi \equiv_{\cal R} \varphi$ (which is in
$\Pi_1^{\sf P}$).
However, the $\equiv_{\cal D}$-decision problem of ${\it Inv}_{3,1}$ is 
$\oplus {\sf P}$-complete, hence  $\oplus {\sf P}$ would be contained in
$\Sigma_2^{\sf P}$. By \cite{Toda91},
$\oplus {\sf P} \subseteq \Sigma_2^{\sf P}$ implies that {\sf PH} collapses
to $\Pi_3^{\sf P} \cap \Sigma_3^{\sf P}$.
 \ \ \ $\Box$

\medskip

The remarks after the proof of Theorem \ref{pivotLength_GammaTau} 
apply also to ${\it Inv}_{k,1}$ when $k \neq 3$.

\bigskip

Finally, we consider the generalized word problem of ${\it Inv}_{k,1}$ or 
$G_{k,1}$ in $M_{k,1}$ over a finite generating set $\Gamma_M$ of $M_{k,1}$ 
(or over a circuit-like generating set $\Gamma_M \cup \tau$).
These generalized word problems over $\Gamma_M$ or over 
$\Gamma_M \cup \tau$, are specified as follows. \\  
{\sf Input:} \ $\varphi \in M_{k,1}$, given by a word over $\Gamma_M$ (or 
over $\Gamma_M \cup \tau$). \\  
{\sf Question (generalized word problem of ${\it Inv}_{k,1}$ in 
$M_{k,1}$):} \ Is $\varphi$ in ${\it Inv}_{k,1}$ ?  \\   
{\sf Question (generalized word problem of $G_{k,1}$ in $M_{k,1}$):} 
 \ Is $\varphi$ in $G_{k,1}$ ?

\begin{pro} \label{genwpInvGinM_fingen} \ When inputs are given over a finite 
generating set of $M_{k,1}$ the generalized word problems of 
${\it Inv}_{k,1}$ and of $G_{k,1}$ in $M_{k,1}$ are in {\sf P}.
\end{pro}
{\bf Proof.} \ We saw in \cite{BiRL}, Sections 5.1 and 6.1, that for 
inputs over $\Gamma$ we can check $\equiv_{\cal R}$ and $\equiv_{\cal L}$ 
in deterministic polynomial time. Since $\varphi \in G_{k,1}$ iff 
$\varphi \equiv_{\cal L} {\bf 1}$ and $\varphi \equiv_{\cal R} {\bf 1}$,
it follows that the generalized word problem of $G_{k,1}$ in $M_{k,1}$ 
over $\Gamma$ is in {\sf P}.

To solve the generalized word problem of ${\it Inv}_{k,1}$ in $M_{k,1}$, 
observe that $\varphi \in {\it Inv}_{k,1}$ iff $\varphi$ is injective,
which holds iff for evey $y \in {\sf imC}(\varphi)$, 
$|\varphi^{-1}(y)| = 1$.
Recall that when $\varphi$ is given by a word over $\Gamma$, we can 
compute ${\sf imC}(\varphi)$ as an explicit list of words, in 
deterministic polynomial time (Corollary 4.11 in \cite{BiThomMon}).
Also, for each $y \in {\sf imC}(\varphi)$ we can compute a finite-state 
automaton ${\cal A}_y$ accepting $\varphi^{-1}(y)$; the set 
$\varphi^{-1}(y)$ is finite and ${\cal A}_y$ is acyclic, and reduced. 
We have $|\varphi^{-1}(y)| = 1$ iff every state in ${\cal A}_y$ has 
out-degree 1, i.e., the graph of ${\cal A}_y$ is a chain (with every
edge labeled by one letter). Checking whether ${\cal A}_y$ is a chain
can be done in polynomial time, so the generalized word problem of 
${\it Inv}_{k,1}$ in $M_{k,1}$ is in {\sf P}. 
 \ \ \ $\Box$

\begin{pro} \label{genwpInvinM_infgen} \ When inputs are given over a 
circuit-like generating set $\Gamma \cup \tau$ of $M_{k,1}$ the 
generalized word problem of ${\it Inv}_{k,1}$ in $M_{k,1}$ is 
{\sf coNP}-complete.
\end{pro}
{\bf Proof.} \ An element $\varphi \in M_{k,1}$ belongs to 
${\it Inv}_{k,1}$ iff $\varphi$ is injective. 
By Theorem 4.5 in \cite{BiThomMon}, the length of the longest words in 
${\sf domC}(\varphi)$ is $\leq c \cdot |\varphi|_{\Gamma \cup \tau}$,
for some constant $c$. Hence, non-injectiveness of $\varphi$ can be 
decided in nondeterministic polynomial time by guessing two different
words $x_1, x_2 \in {\sf domC}(\varphi)$ and checking that 
$\varphi(x_1) = \varphi(x_2)$. We know from Theorem 4.12 in 
\cite{BiThomMon} that $\varphi(x_1)$ and $\varphi(x_2)$ can be computed 
in deterministic polynomial time.
Hence, the generalized word problem of ${\it Inv}_{k,1}$ in $M_{k,1}$
is in {\sf coNP}.  

In Prop.\ 6.5 in \cite{BiRL} it was proved that the injectiveness 
problem for combinational circuits is {\sf coNP}-complete.
Since combinational circuits can be represented by words over the 
circuit-like generating set $\Gamma_M \cup \tau$ of $M_{k,1}$, it follows 
that the generalized word problem of ${\it Inv}_{k,1}$ in $M_{k,1}$ is 
{\sf coNP}-hard. 
 \ \ \ $\Box$

\bigskip

\noindent {\bf Open question:} \ What is the complexity of the 
generalized word problem of the Thompson-Higman group $G_{k,1}$ in
$M_{k,1}$ when the input is a word over a circuit-like generating set 
$\Gamma_M \cup \tau$ of $M_{k,1}$?

\smallskip

We know that the problem is in $\Pi_2^{\sf P}$, since the question
whether $\varphi \equiv_{\cal L} {\bf 1}$ is in {\sf coNP} (by Theorem 6.7 
in \cite{BiRL}), and the question whether 
$\varphi \equiv_{\cal R} {\bf 1}$ is in $\Pi_2^{\sf P}$. 

We also know that the problem is {\sf coNP}-hard.
Indeed, Prop.\ 6.5 in \cite{BiRL} gives a polynomial-time reduction 
$B \mapsto F_B$ where $B$ is any boolean formula, and $F_B \in M_{k,1}$ 
(given by a word over $\Gamma_M \cup \tau$) is such that:
  
(1) \ if $B$ is a tautology then $F_B = {\bf 1}$ as an element of
   $M_{k,1}$;

(2) \ if $B$ is not a tautology then $F_B$ is not injective. 

\noindent Since ${\bf 1} \in G_{k,1}$, whereas non-injective elements 
are not in $G_{k,1}$, this reduces the tautology problem (which is 
{\sf coNP}-complete) to the generalized word problem of $G_{k,1}$ in 
$M_{k,1}$. 

\bigskip

\noindent {\bf Open question:} \ What is the distortion of $G_{k,1}$ 
(over $\Gamma_G \cup \tau$) within $M_{k,1}$ (over $\Gamma_M \cup \tau$)?  
Similarly, what is the distortion of ${\it Inv}_{k,1}$ in $M_{k,1}$
over circuit-like generating sets?

%The distortion of $G_{k,1}$ (or ${\it Inv}_{k,1}$) in $M_{k,1}$ over
%finite generating sets has a polynomial upper-bound.

%%%%%%%%%%%%%%%%%%%%%%%%%%%%%%%%%%%%%%%%%%%%%%%%%%%%%%%%%%%%%%%%%%%%%%
%% Section: Appendix
%%%%%%%%%%%%%%%%%%%%%%%%%%%%%%%%%%%%%%%%%%%%%%%%%%%%%%%%%%%%%%%%%%%%%%

\section{Appendix: Search problems, {\sf NPsearch}, and {\sf xNPsearch}}

A {\it search problem} is a relation of the form
$R \subseteq A^* \times B^*$ (where $A$ and $B$ are finite alphabets).
We usually formulate the search problem $R$ in the following form. \\
{\sf Input:} \ A string $x \in A^*$; \\
{\sf Premise:} \ There exists $y \in B^*$ such that $(x,y) \in R$; \\
{\sf Search:} \ Find one $y \in B^*$ such that $(x,y) \in R$.

\smallskip

A {\it premise problem} is a decision problem (or a search problem) which,
in addition to an input and a question (or a requested output), also has
a premise concerning the input.
The premise, also called ``pre-condition'', is an assumption about the
input that any algorithm for the problem is allowed to use as a fact.
The algorithm does not need to check whether the assumption actually holds
for the given input and, indeed, {\it we don't care} about the answer (or 
the output) when the premise does not hold. In the literature the word 
``promise'' is often used for ``premise'' (although, according to the 
dictionaries of the English language, ``premise'' is more logical).

\smallskip

For $R \subseteq A^* \times B^*$, the domain of $R$ is
 \ ${\sf Dom}(R) \ = \ $
$\{x \in A^* : (\exists y \in B^*) [(x,y) \in R]\}$; in words,
${\sf Dom}(R)$ is the set of inputs for which the search problem $R$ has
a solution.
The membership problem of ${\sf Dom}(R)$ is called the {\em decision
problem} associated with the search problem $R$.
The membership problem of $R$ is called the {\em verification problem}
associated with $R$.

\smallskip

The best known search problem complexity class is {\sf NPsearch} (called 
{\sf FNP} or
``function problems associated with {\sf NP}'' in \cite{Papadim}).
The class {\sf NPsearch} consists of all relations of the form
$R \subseteq A^* \times B^*$ (where $A$ and $B$ are finite alphabets)
such that (according to \cite{Papadim}, pages 227-240):  \\
(1) \ \ \ the membership problem of $R$ (i.e., the verification problem)
belongs to {\sf P}, i.e., there exists a deterministic polynomial-time
algorithm which on input $(x,y) \in A^* \times B^*$ decides whether
$(x,y) \in R$; \\
(2) \ \ \ $R$ is polynomially balanced; this means that there exists a 
polynomial $p$ such that for all $(x,y) \in R$: \ $|y| \leq p(|x|)$.

\smallskip

When $R$ is in {\sf NPsearch} then the associated decision problem is in
{\sf NP}. The complementary decision problem (namely the task of
answering ``no'' on input $x$ iff there is no $y$ such that $(x,y) \in R$)
is in {\sf coNP}. It is interesting to
compare the class {\sf NPsearch} also with $\#{\sf P}$, which consists of
the functions that count the number of solutions of {\sf NPsearch}
problems.

By definition, a deterministic algorithm $\cal A$ solves the search
problem $R$ iff for every input $x \in {\sf Dom}(R)$, the algorithm
$\cal A$ outputs an element $y \in B^*$ such that $(x,y) \in R$.
No requirement is imposed on $\cal A$ when $x \not\in {\sf Dom}(R)$;
however, if complexity bounds are known (or required) for $\cal A$ the
above definition implies that $\cal A$ also, indirectly, determines
whether $x \not\in {\sf Dom}(R)$, and we output ``no'' in that case.
Probabilistic solutions of a search problem can also be defined. One way
to do that is to say that a probabilistic algorithm $\cal A$ solves the
search problem $R$ \ iff \ for every $x \in {\sf Dom}(R) :$
 \ $P(\{ y \in B^* : y = {\cal A} \ {\rm and} \ (x,y) \in R\})$
 $ \ \geq \ c$ \ (where $c$ is a constant, $0 < c < 1$).
No requirement is imposed on $\cal A$ when $x \not\in {\sf Dom}(R)$. But
since $R$ is in {\sf P}, proposed false solutions can be ruled out.

\smallskip

\noindent {\bf Remark.} \ The idea of solving a search problem by a
deterministic algorithm explains why {\sf NPsearch} was called {\sf FNP}
(where the ``{\sf F}'' stands for ``function''). However, it is better not
to attach the word ``function'' to {\sf NPsearch} because the problems in
{\sf NPsearch} are relations.
Functions may play a role in special ways of solving a search problem; but
other, non-functional, solutions of search problems are often considered
too, e.g., probabilistic algorithms.

\smallskip

Following \cite{Papadim}, page 229, we define the concept of a
{\bf polynomial-time many-to-one search reduction} from a search problem
$R_1 \subseteq A_1^* \times B_1^*$ to a search problem
$R_2 \subseteq A_2^* \times B_2^*$ as follows. Such a reduction is a 
triple of polynomial-time computable total functions
$\rho_{\rm in}: A_1^* \to A_2^*$,
 \ $\rho_{\rm sol}: A_1^* \times B_2^* \to B_1^*$, and
$\rho_{\rm ver}: A_1^* \times B_1^* \to A_2^* \times B_2^*$ \ such that:

\smallskip

\noindent (1) \ For all $x_1 \in {\sf Dom}(R_1):$
 \ $\rho_{\rm in}(x_1) \in {\sf Dom}(R_2)$.

\smallskip

\noindent (2) \ For all $x_1 \in A_1^*$ and all $y_2 \in B_2^*:$
 \ \ $(\rho_{\rm in}(x_1), \, y_2) \ \in \ R_2$ \ \ implies
 \ \ $(x_1, \, \rho_{\rm sol}(x_1,y_2)) \ \in \ R_1$.

\smallskip

\noindent (3) \ For all $(x_1,y_1) \in A_1^* \times B_1^* : $
 \ $(x_1,y_1) \in R_1$ \ iff \ $\rho_{\rm ver}(x_1,y_1) \in R_2$. Moreover,
 $\rho_{\rm ver}$ is ``polynomially balanced'', i.e.,
 there is a polynomial $p$ such that for all
 $(x_1,y_1) \in A_1^* \times B_1^* :$
 \ if $\rho_{\rm ver}(x_1,y_1) = (x_2,y_2)$ then $|y_2| \leq p(|x_2|)$.

\smallskip

\noindent In words, condition (1) says that if $R_1$ has a solution for
input $x_1$ then $R_2$ has a solution for input $\rho_{\rm in}(x_1)$. When
$R_2$ is total, i.e., ${\sf Dom}(R_2) = A_2^*$, then condition (1) holds
automatically.
Condition (2) means that every $R_2$-solution $y_2$ for input
$\rho_{\rm in}(x_1)$ yields an $R_1$-solution $\rho_{\rm sol}(x_1,y_2)$ for
input $x_1$. Condition (3) means that the verification problem of $R_1$
reduces to the verification problem of $R_2$.
As a consequence of condition (3), the class {\sf NPsearch} is closed under
polynomial-time many-to-one search reduction.
(Condition (3) is usually omitted in the literature; however, the literature
also claims that {\sf NPsearch} is closed under search reduction, but this 
does not follow from conditions (1) and (2) alone.)
The pair of maps $(\rho_{\rm in}, \rho_{\rm sol})$ is called the
{\it input-output reduction}, and the map $\rho_{\rm ver}$ is called the
{\it verification reduction}.

\smallskip

There are well-known search problems that are closely related to {\sf NP} 
but that don't exactly fit into the class {\sf NPsearch}.
For example, the search version of integer linear programming is not
polynomially balanced: for some inputs there are infinitely many solutions,
of unbounded size, although there also exist polynomially bounded
solutions for every input that has a solution. Similar examples are
certain versions of the Traveling Salesman problem, or finding solutions
to certain equations (search version of problems on pp.\ 249-253
in \cite{GaJo}).
We prove in Section 5.2 that the special multiplier search problem for 
$\equiv_{\cal J} {\bf 1}$ in $M_{k,1}$ (over $\Gamma \cup \tau$) is another 
example. In the $\equiv_{\cal J} {\bf 1}$ multiplier search problem,
when there are solutions then there are also solutions that are
polynomially bounded and verifiable in deterministic polynomial time. But
the general verification problem for $\equiv_{\cal J} {\bf 1}$ is not 
polynomially balanced.   The main observation is that 
in a search problem we only want to find {\it one} solution for 
each input, so the difficulty of the general verification problem and the
size of all solutions in general should not concern us.

Therefore we introduce the class {\sf xNPsearch} {\bf (extended NP search)},
consisting of all relations of the form $R \subseteq A^* \times B^*$
(where $A$ and $B$ are finite alphabets) such that there is a relation
$R_0 \subseteq R$ with the properties \\
(1) \ \ \ $R_0 \in {\sf NPsearch}$, \\
(2) \ \ \ ${\sf Dom}(R) = {\sf Dom}(R_0)$.

\smallskip

When $R$ is in {\sf xNPsearch} then the associated decision problem is in
{\sf NP}, just as for {\sf NPsearch}, since problems in {\sf NPsearch} and
{\sf xNPsearch} have the same domains.

By definition, a search problem $R$ is {\sf NPsearch}{\bf -complete} iff
$R$ is in {\sf NPsearch}, and every problem in {\sf NPsearch} can be reduced
to $R$ by a polynomial-time many-to-one search reduction. We say that $R$ 
is {\sf xNPsearch}-complete iff there is an {\sf NPsearch}-complete 
problem $R_0$ such that $R_0 \subseteq R$ and 
${\sf Dom}(R_0) = {\sf Dom}(R)$.

It follows that an {\sf xNPsearch}-complete problem is in
{\sf xNPsearch}.
And it follows that an {\sf NPsearch}-complete problem is automatically
{\sf xNPsearch}-complete.

\smallskip

An example of an {\sf NPsearch}-complete (hence {\sf xNPsearch}-complete)
problem is the following, called {\sc SatSearch}; it is the relation
 \ $\{ (B,t) : B $ is a boolean formula with $m$ variables, \ $m > 0$,
 \ $t \in \{0,1\}^m$, and \ $B(t) = 1\}$.
Equivalently, {\sc SatSearch} is specified as follows: \\
{\sf Input:} \ A boolean formula $B(x_1, \ldots, x_m)$ (where $m$ is part
of the input, hence variable). \\
{\sf Premise:} \ $B(x_1, \ldots, x_m)$ is satisfiable. \\
{\sf Search:} \ Find a satisfying truth-value assignment $t \in \{0,1\}^m$
for $B(x_1, \ldots, x_m)$.

%%%%%%%%%%%%%%%%%%%%%%%%%%%

\bigskip

\bigskip

\noindent {\bf Acknowledgements:}  I would like to thank Lane Hemaspaandra 
for references on counting complexity classes, and Sunil Shende, for 
clarifications on reductions between search problems.

%%%%%%%%%%%%%%%%%%%%%%%%%%%%%%%%%%%%%%%%%%%%%%%%%%%%%%%%

\bigskip

\bigskip

%%%%%%%%%%%%%%%%%%%%%%%%%%%%%%%%%%%%%%%%%%%%%%%%%%%%%%%%%%%%%%%%%%%%%%%%%%%%%%%
%% \small

%%%%%%%%%%%%%%%%%%%%%%%%%%%%%

\bigskip

\bigskip

\noindent {\bf Jean-Camille Birget} \\
Dept.\ of Computer Science \\
Rutgers University at Camden \\
Camden, NJ 08102, USA \\
{\tt birget@camden.rutgers.edu}


{\small 

\begin{thebibliography}{99}

\bibitem{BeiGill} R.\ Beigel, J.\ Gill, ``Counting classes: Thresholds,
  parity, mods, and fewness'', {\it Theoretical Computer Science} 103(1)
  (1992) 2-23.

\bibitem{BeiGillHert} R.\ Beigel, J.\ Gill, U.\ Hertrampf, ``Counting
  classes: Thresholds, parity, mods, and fewness'', {\it Proc.\ 7th
  Annual Symposium on Theoretical Aspects of Computer Science} (STACS'90)
  LNCS 415 Springer Verlag (1990), pp.\ 49-57.

\bibitem{BiRL} J.C.Birget, ``The $\cal R$- and $\cal L$-orders of the
  Thompson-Higman monoid $M_{k,1}$ and their complexity'', Mathematics ArXiv,
  http://arXiv.org/abs/0812.4434 [math.GR], 23 Dec 2008.

\bibitem{BiThomMon} J.C.~Birget, ``Monoid generalizations of the Richard
  Thompson groups'', {\it J.~of Pure and Applied Algebra}, 213(2) (Feb.\
  2009) 264-278. (Online pre-publication DOI:
  http://dx.doi.org/10.1016/j.jpaa.2008.06.012 )
  (Preprint: \ Mathematics ArXiv
  http://arXiv.org/abs/math.GR/0704.0189 , April 2007.)

\bibitem{BiDistor} J.C.\ Birget, ``One-way permutations, computational
  asymmetry and distortion'', {\it J.~of Algebra}, 320(11) (Dec.\ 2008)
  4030-4062. \ (Online pre-publication DOI:
  http://dx.doi.org/10.1016/j.jalgebra.2008.05.035 )
  (Preprint: \ Mathematics ArXiv
  http://arxiv.org/abs/0704.1569,  April 2007).

\bibitem{BiFact} J.C.\ Birget, ``Factorizations of the Thompson-Higman
  groups, and circuit complexity'', {\it International J.~of Algebra and
  Computation}, 18.2 (March 2008) 285-320.
  (Preprint: \  Mathematics ArXiv
  http://arXiv.org/abs/math.GR/0607349, July 2006.)

\bibitem{BiCoNP} J.C.~Birget, ``Circuits, coNP-completeness, and the groups
  of Richard Thompson'', {\it International J.~of Algebra and Computation}
  16(1) (Feb.\ 2006) 35-90. \\  
  (Preprint: 
  Mathematics ArXiv http://arXiv.org/abs/math.GR/0310335, Oct.\ 2003).

\bibitem{BiThomps} J.C.\ Birget, ``The groups of Richard Thompson and
  complexity'', {\it International J. of Algebra and Computation} 14(5,6)
  (Dec.\ 2004) 569-626.
  (Preprint: \  Mathematics ArXiv math.GR/0204292, Apr.\ 2002).

\bibitem{CGHHSWW} J.-Y.\ Cai, T.\ Gundermann, J.\ Hartmanis,
  L.\ Hemachandra, V.\  Sewelson, K.\ Wagner, G.\ Wechsung,
  ``The boolean hierarchy I: structural properties'',
  {\it SIAM Journal on Computing} 17(6) (1988) 1232-1252;
  ``The boolean hierarchy II: applications'',
  {\it SIAM Journal on Computing} 18(1) (1089) 95-111.  

\bibitem{CaiHem90} J.-Y.\ Cai, L.\ Hemachandra, ``On the power of parity
  polynomial time'', {\it Mathematical Systems Theory} 23(2) (1990) 95-106.

\bibitem{ChangKadin} R.\ Chang, J.\ Kadin, ``The boolean hierarchy and
  the polynomial hierarchy: a closer connection'', {\it SIAM Journal on
  Computing}, 25(2) (1996) 340-354.

\bibitem{CliffPres} A.H.\ Clifford, G.B.\ Preston, {\it The algebraic
  theory of semigroups}, Vol.\ 1 (Mathematical Survey, No 7 (I))  American
  Mathematical Society, Providence (1961).

\bibitem{DuKo} D.Z.\ Du, K.I.\ Ko, {\it Theory of computational complexity},
  Wiley (2000).

\bibitem{DurandHerKol} A.\ Durand, M.\ Hermann, Ph.\ Kolaitis,
  ``Subtractive reductions and complete problems for counting classes'',
  {\it Theoretical Computer Science} 340 (2005) 496-513.

\bibitem{GaJo} M.\ Garey, D.\ Johnson, {\it Computers and intractability,  
   a guide to the theory of NP-completeness}, Freeman (1979).

\bibitem{GP86} L.\ Goldschlager, I.\ Parberry, ``On the construction of
  parallel computers from various bases of boolean functions'',
  {\it Theoretical Computer Science} 43 (1986) 43-58.

\bibitem{Grillet} P.A.\ Grillet, {\it Semigroups, An introduction to the
  structure theory}, Marcel Dekker, New York (1995).

\bibitem{HemOgiCompan} L.\ Hemaspaandra, M.\ Ogihara, {\it The complexity
  theory companion}, Springer (2002).

\bibitem{HemaspVollmer} L.\ Hemaspaandra, H.\ Vollmer, ``The satanic
   notations: Counting classes beyond \#P and other definitional
   adventures'', {\it SIGACT News} 26(1) (1995) 2-13.

\bibitem{Hig74} G.\ Higman, ``Finitely presented infinite simple groups'',
  Notes on Pure Mathematics 8, The Australian National University,
  Canberra (1974).

\bibitem{Kadin} J.\ Kadin, ``The polynomial hierarchy collapses if the
  boolean hierarchy collapses'', {\it SIAM J.\ on Computing} 17 (1988)
  1263 - 1282. (Corrections in \cite{ChangKadin}.)

\bibitem{KoSchoTo} J.\ K\"obler, U.\ Sch\"oning, J.\ Tor\'an, ``On
  counting and approximation'', {\it Acta Informatica} 26(4) (1989) 363-379.

\bibitem{Lallement} G.\ Lallement, {\it Semigroups and combinatorial 
  applications}, Wiley (1979).

\bibitem{Papadim} Ch.\ Papadimitriou, {\it Computational complexity},
  Addison-Wesley (1994).

\bibitem{PapadYannak} Ch.\ Papadimitriou, M.\ Yannakakis, ``The complexity
  of facets (and some facets of complexity), {\it Journal of Computer
  and System Sciences} 28 (1984) 244-259. [Conference version: ACM STOC
  (1982) 255-259.]

\bibitem{PapadZachos} Ch.\ Papadimitriou, S.\ Zachos, ``Two remarks on the
  power of counting'', {\it Proc.\ 6th GI Conf. on Theoretical Computer
  Science}, Springer-Verlag LNCS No.\ 149 (1983) pp.\ 269-276.

\bibitem{Scott} Elizabeth A. Scott, ``A construction which can be used
  to produce finitely presented infinite simple groups'',
  {\it J. of Algebra} 90 (1984) 294-322.

\bibitem{ScottTour} Elizabeth A. Scott, ``A tour around finitely presented 
  infinite simple groups'', in {\it Algorithms and classification in
  combinatorial group theory}, G.\ Baumslag and C.F.\ Miller III editors,
  MSRI Publ.\ vol.\ 23, Springer Verlag (1992)

\bibitem{TodaDiss} S.\ Toda, ``Computational complexity of counting
  complexity classes'', PhD thesis, Tokyo Institute of Technology, 
  Tokyo, Japan (1991).

\bibitem{Toda91} S.\ Toda, ``PP is as hard as the polynomial hierarchy'', 
  {\it SIAM J.\ on Computing} 20 (1991) 865-877.

\bibitem{Valiant1} L.G.\ Valiant, ``The complexity of computing the
  permanent'', {\it Theoretical Computer Science}, 8(2) (1979) 189-201.

\bibitem{Valiant2} L.G.\ Valiant, ``The complexity of enumeration and
  reliability problems'', {\it SIAM J.\ on Computing}, 3(3) (1979) 410-421.

\bibitem{Handb} J.\ van Leeuwen (editor), {\it Handbook of theoretical
  computer science}, volume {\bf A}, MIT Press and Elsevier (1990).


\end{thebibliography}

}
\end{document}